\documentclass[letterpaper, 10 pt, journal, twoside]{IEEEtran}

\IEEEoverridecommandlockouts                            

\usepackage[pdftex]{graphicx}
\graphicspath{{figure/}}
\usepackage{amsmath}
\usepackage{amssymb}
\usepackage{url}

\usepackage{amsthm} 
\usepackage{cite}
\usepackage{color}
\usepackage{comment}
\usepackage{fancyhdr} 

\usepackage{algorithm}
\usepackage{algorithmic}
\usepackage{bm}
\usepackage{bbm}

\newtheorem{theorem}{\textbf{Theorem}}
\newtheorem{lemma}{\textbf{Lemma}}
\newtheorem{assumption}{\textbf{Assumption}}
\newtheorem{remark}{\textbf{Remark}}

\newcommand{\bx}{\bm{x}}

\newcommand{\by}{\bm{y}}
\newcommand{\bz}{\bm{z}}
\newcommand{\bp}{\bm{p}}
\newcommand{\bs}{\bm{s}}
\newcommand{\bg}{\bm{g}}
\newcommand{\bv}{\bm{v}}

\newcommand{\bd}{\bm{d}}

\newcommand{\bh}{\bm{h}}
\newcommand{\bla}{\bm{\lambda}}
\newcommand{\bom}{\bm{\omega}}
\newcommand{\R}{\mathbb{R}}

\newcommand{\sX}{\mathcal{X}}

\newcommand{\sA}{\mathcal{A}}
\newcommand{\sO}{\mathcal{O}}
\newcommand{\sZ}{\mathcal{Z}}

\newcommand{\sF}{\mathcal{F}}
\newcommand{\sC}{\mathcal{C}}
\newcommand{\sM}{\mathcal{M}}

\newcommand{\B}{\mathbb{B}}

\newcommand{\bxi}{\bm{\xi}}

\newcommand{\bmu}{\bm{\mu}}

\newtheorem{corollary}{Corollary}
\newtheorem{proposition}{Proposition}

\newtheorem{definition}{Definition}

\allowdisplaybreaks 

\title{ 
Model-Free Optimal Voltage Control \\ via
Continuous-Time Zeroth-Order Methods }

\author{Xin Chen, Jorge I. Poveda, Na Li
\thanks{X. Chen and N. Li are with the School of Engineering and Applied Sciences, Harvard University, USA; Emails: chen\_xin@g.harvard.edu, nali@seas.harvard.edu.}
\thanks{
J. I. Poveda is with  the Department of Electrical, Computer, and Energy Engineering at the University of Colorado, Boulder, USA; Email: jorge.poveda@colorado.edu.}
\thanks{ 
The work was supported by NSF CNS 1947613, NSF CAREER: ECCS-1553407 and
NSF EAGER: ECCS-1839632.} 
}

\begin{document}

\maketitle

	
\begin{abstract}
 In power distribution systems, the growing penetration of renewable energy resources brings new challenges to maintaining voltage safety, 
 which is further complicated by the limited model information of distribution systems. 
To address these challenges, we develop a model-free optimal voltage control algorithm 
	 based on  projected primal-dual gradient dynamics and  continuous-time zeroth-order method (extreme seeking control). This proposed algorithm i) operates purely based on voltage measurements and does not require any other model information, ii) can drive the voltage magnitudes back to the acceptable range, iii) satisfies the power capacity constraints all the time, iv) minimizes the total operating cost,
	 and v) is implemented in a decentralized fashion where the privacy of controllable devices is preserved and plug-and-play operation is enabled. We prove that the proposed algorithm is semi-globally practically asymptotically stable and is structurally robust to measurement noises.
Lastly, the performance of the proposed algorithm is further  demonstrated via numerical simulations. 

\end{abstract}


\section{Introduction} \label{sec:introduction}

Voltage control in a distribution system aims to maintain the  voltage magnitudes across the power network within an acceptable range \cite{8636257}.  
With  rapidly increasing penetration of renewable energy resources, such as 
 photovoltaic (PV) and  wind generation,  it brings emerging operational challenges to the task of voltage control.  On the one hand,
 the caused reverse power flow  may lead to frequent  over-voltage issues. On the other hand, large-scale renewable generations introduce significant uncertainty and volatility to the distribution systems, making it much harder to model and control.


There have been a large amount of researches \cite{zheng2017robust,7039295,7244261,8779692,8268542,7361761,7028508} devoted to voltage control by regulating the slow time-scale devices (such as voltage regulators, shunt capacitors, and on-load-tap-changer transformers) and fast time-scale devices (such as distributed generations (DGs) and static Var compensators (SVCs)). However, most existing voltage control methods are based on power flow models and assume good knowledge of the distribution systems. Therefore, these methods may not perform well when  such models and information are absent. 
References \cite{8779692,8268542,7361761,7028508} propose feedback
 voltage control schemes based on primal-dual gradient methods, dual ascent approaches, or integral control. Due to the feedback mechanism,
 these schemes circumvent some of the system information, e.g., 
  real-time uncontrollable power injections, while the distribution network model, such as line parameters and network topology, is still required.  In practice,  high-accuracy network models and onsite identified network parameters are unavailable for many   distribution systems. 
Moreover, network reconfiguration, line faults, and other operational factors also change the system model from time to time. Hence, it is desirable for the voltage control schemes to  operate well in the absence of  system models and adapt fast to  time-varying operational conditions.

The  deployment of smart meters and upgraded communication infrastructures 
offer  an   opportunity  to  overcome  these  challenges through  real-time   monitoring  and  control, which motivates the data-driven voltage control techniques. A  type of such data-driven schemes \cite{8412143,7348696,7540852} is to approximate the nonlinear power flow relation with a linear sensitivity model (e.g., the LinDistflow model \cite{7361761}), and then 
online estimate  the  model using  measurements and regression methods for voltage control.  These schemes generally require a control center to  store a large amount of measurement data and 
solve high-dimensional regression problems in real-time. Reference  \cite{8873667} proposes to 
reduce the complexity of the linear regression by assuming and exploiting the knowledge of network topology and line resistance-to-reactance ratios. 
The other type of data-driven schemes is the end-to-end \emph{model-free} control, such as reinforcement learning (RL), which does not explicitly estimate the system model and 
makes  decisions directly based on measurements.
 A number of recent works \cite{8944292,9143169,9076841,9274529}  propose to learn voltage control policies using various RL techniques;  see review article \cite{chen2021reinforcement} and references therein for a more comprehensive view.  However, applying RL to the control of physical systems is still under development and generally has many  limitations, such as safety problems (e.g., physical constraint violation), scalability issues, unstable training process, limited or no theoretical guarantee, etc.

An alternative type of model-free control is based on zeroth-order (or gradient-free) methods \cite{9147814}. 
In particular, \emph{extremum seeking} (ES) control \cite{ariyur2003real} is a classic continuous-time zeroth-order optimal control method, which operates using only  the output measurements.
ES control  attracts surging recent attention and has been
 applied in   broad power system applications, including
 energy consumption control \cite{7397989}, voltage phasor regulation \cite{8888186},  maximum power point tracking \cite{6362193}, etc.  Moreover,
references \cite{7350258, 8772122}  develop ES control algorithms to modulate the power injections of distributed energy resources   for voltage regulation.  In \cite{8478426},  hardware-in-the-loop experiments are conducted to verify the viability of a  ES-based voltage control scheme. 
Despite  these progresses,
one major limitation of existing ES algorithms is that the  constraints are not well addressed. Most of the ES methods above 
consider   unconstrained 
 optimization problems for simplicity or penalize the constraint violation in the objective. However,
 there are various physical constraints, e.g., the power capacity limits, that need to be enforced in practice.

\textbf{Contributions.} 
In this paper, we study the real-time voltage control through modulating the active and reactive power outputs of fast time-scale controllable devices.  
To overcome the  challenges described above,
we develop a  model-free optimal voltage control  algorithm  based on  projected primal-dual gradient dynamics (P-PDGD) and ES control. Specifically, 
by leveraging the structure of P-PDGD, the proposed algorithm can  steer the system to an optimal operating point while satisfying the physical constraints. Then  ES control is adopted to make this algorithm ``model-free" in the sense that the distribution system model is circumvented. The main merits of the proposed algorithm are explained as below:
\begin{itemize}
    \item [1)] (\emph{Optimality}). The proposed algorithm can  drive the voltage magnitudes back to the acceptable range while minimizing the total operating cost and always satisfying the power  capacity  constraints. 
    \item[2)] (\emph{Model-Free}). The proposed algorithm is an end-to-end model-free control method that operates purely based on the voltage measurements from the monitored buses. The model information of  distribution networks and other power injections  is not needed.
    \item[3)] (\emph{Adaptive}). By exploiting real-time measurement, this algorithm is a feedback mechanism that can adapt fast to changes in the dynamical system environment. 
    \item[4)] (\emph{Decentralized}). This  algorithm is implemented in a decentralized manner, where the privacy of each device can be preserved. Moreover, it allows plug-and-play operation and thus is robust to single/multi-point failures.
    \item[5)] (\emph{Guaranteed Performance}). We mathematically prove the
 semi-global   practical   asymptotical   stability and  the structural  robustness  (to small measurement  noise) of the proposed algorithm, and  numerically verify its effectiveness, optimality and robustness via  simulations.
\end{itemize}

To the best of our knowledge, this is the first work on voltage control that  unifies all the above features. We also emphasize that 
 the proposed ES-P-PDGD algorithm is a generic model-free method that can be 
 applied to many other multi-agent optimization and control problems. Comparing with existing  ES methods, our algorithm can enforce hard physical constraints without sacrificing other performance.

Lastly, we mention a closely-related work \cite{9147814}. It proposes a model-free primal-dual projected gradient algorithm for real-time optimal power flow  based on discrete-time zeroth-order methods, 
but it makes relatively strong assumptions on the
problem setting and lacks explicit convergence results. In
contrast, this paper uses and studies the continuous-time ES
control dynamics, and provides clear stability guarantee. Besides, distinguished from the projection method used in \cite{9147814}, our algorithm employs the global projection and it leads to  a \emph{Lipschitz continuous} projected dynamical system (see Remark \ref{remark:project}), which facilitates  the theoretical analysis.


The remainder of this paper is organized as follows: Section \ref{sec:problem} presents the optimal voltage control problem and the preliminaries on ES control. Section \ref{sec:algorithm} develops the model-free algorithm based on P-PDGD and ES. 
 Section \ref{sec:performance} analyzes the theoretical performance of the proposed algorithm.
  Numerical tests are conducted in Section \ref{sec:simulation}, and conclusions are drawn in Section \ref{sec:conclusion}.

\textbf{Notations.} 
We use unbolded lower-case letters for scalars, and bolded lower-case letters for column vectors. 
$\mathbb{R}_+:=[0,+\infty)$ denotes the set of non-negative  real values. $|\cdot|$  denotes the cardinality of a set. $||\cdot||$ denotes the L2-norm of a vector. $[\bx; \by] := [\bx^\top, \by^\top]^\top$ denotes the column merge of  vectors $\bx,\by$. $\blacktriangle$  highlights the  definition of new notations.

\section{Problem Formulation and  Preliminaries} \label{sec:problem}

In this section, we present the formulation of the optimal voltage control problem and introduce the preliminaries on extremum seeking control. 

\subsection{Optimal Voltage Control Problem}

Consider a distribution network with the monitored bus set $\sM$ and the controllable device set $\sC$. Each bus $j\in\sM$ has real-time voltage measurement, and the power injection of each device
$i\in\sC$ can be adjusted for voltage regulation.
Depending on the practical system configuration, the controllable devices are flexible to locate at any buses of the distribution network.
The optimal voltage control (OVC) problem is formulated as model (\ref{eq:ovc}) and explained below:
\begin{subequations}   \label{eq:ovc}
\begin{align}
    \text{Obj.}\ & \min_{\bx} \sum_{i\in\sC} c_i(\bx_i)  \label{eq:ovc:obj} \\
\text{s.t.}\ &  
 \ \bx_i \in \mathcal{X}_i, && i\,\in\, \mathcal{C} \label{eq:ovc:x}\\
& \      \underline{v}_j\leq v_j(\bx) \leq \bar{v}_j, &&j\in \mathcal{M}. \label{eq:ovc:v}
\end{align}
\end{subequations}

\subsubsection{Decision Variable and its Feasible Set}

The decision variable $\bx_i$ is the power injection of controllable device $i\in\sC$, and its power capacity constraints are described with the feasible set $\sX_i$ in \eqref{eq:ovc:x}. We define 
\begin{align*}
    \bx:=(\bx_i)_{i\in\sC},\quad \sX:= \prod_{i\in\sC} \sX_i\, .
\end{align*}
Specifically, we consider the following two types of  devices for real-time voltage control with $\sC = \sC_{\mathrm{svc}}\cup \sC_{\mathrm{dg}}$:

i) \emph{Static Var Compensator (SVC)} with the reactive power injection $\bx_i:= q_i$ and the power capacity constraint (\ref{eq:svc:cap}):
\begin{align}\label{eq:svc:cap}
     \mathcal{X}_i:= \{\bx_i|\, \underline{q}_i\leq q_i\leq \bar{q}_i  \}, \quad i\in \mathcal{C}_{\mathrm{svc}}
\end{align}
where $ \bar{q}_i$ and $\underline{q}_i$ are the upper and lower limits, respectively.

ii) \emph{Distributed Generation (DG)} with the active and reactive power injection $\bx_i:= [p_i,q_i]^\top$  and  constraint (\ref{eq:dg:cap}):
\begin{align}\label{eq:dg:cap}
   \mathcal{X}_i:= \{\bx_i|\, \underline{p}_i\leq p_i\leq \bar{p}_i, p_i^2+q_i^2\leq \bar{s}_i^2 \},\quad i\in \mathcal{C}_{\mathrm{dg}}
\end{align}          
where  $\bar{p}_i$ and $\underline{p}_i$ are the upper and lower limits of active power, and $\bar{s}_i$ denotes the apparent power capacity.

\subsubsection{Network and Voltage Constraints} 

$v_j$ in (\ref{eq:ovc:v}) denotes the voltage magnitude  at bus $j\in\sM$, and $\underline{v}_j$ and $\bar{v}_j$ are the lower and upper voltage limits, respectively.
We use the functional form $v_j(\bx)$ to describe
 the input-output map from the controllable power injection $\bx$ to the voltage magnitude $v_j$.  Essentially,  $\bv(\bx)\!:=\!(v_j(\bx))_{j\in\sM}$
captures the nonconvex power flow relation,  distribution network model, and other uncontrollable power injections; see \cite{8779692, 8873667} for details.

By ``system model", we specifically refer to  function $\bv(\bx)$. And ``model-free" means 
that the formulation of  ${\bv}(\bx)$ is unknown and no model estimation is performed for it.


\subsubsection{Objective Function}

The objective (\ref{eq:ovc:obj}) aims to minimize the total operating cost with the cost function $c_i(\cdot)$ for each device $i\in\mathcal{C}$. For instance, the quadratic function \eqref{eq:cost} is a widely used objective \cite{8873667,8779692}:
\begin{align}\label{eq:cost}
    c_i(\bx_i) =\begin{cases}
    c_i^{\mathrm{svc}}\cdot q_i^2, & i\in \mathcal{C}_{\mathrm{svc}}\\
    c_{p,i}^{\mathrm{dg}}\cdot p_i^2 + c_{q,i}^{\mathrm{dg}}\cdot q_i^2, & i\in \mathcal{C}_{\mathrm{dg}}
    \end{cases}
\end{align}
where $c_i^{\mathrm{svc}}, c_{p,i}^{\mathrm{dg}}, c_{q,i}^{\mathrm{dg}}$ are the cost coefficients.

 We summarize the known and unknown information in our problem setting with the following assumption.
\begin{assumption}\label{ass:known}
The gradient of individual cost function, i.e., $\nabla c_i(\cdot)$, exists and is known to each device $i\in\sC$ itself, as well as 
 the feasible set $\sX_i$. The function $\bv(\bx)$ is unknown but the real-time measurement of $\bv$ is available.
\end{assumption}

\begin{remark} 
\normalfont In the follows, we consider a general convex cost function $c_i(\cdot)$, while the quadratic cost function (\ref{eq:cost})  is only adopted for simulations. Besides,  we 
assume that the gradient  $\nabla c_i(\cdot)$  is known to each device   for simplicity.
Nevertheless, the proposed voltage control algorithm is applicable to the case when the gradient $\nabla c_i(\cdot)$  is unknown but the cost value $c_i$ can be measured in real time. Similarly, if the real-time measurement  of network loss is available, the cost of network loss $l(\bx)$ can be also included in  objective (\ref{eq:ovc:obj}).
\end{remark}

\subsection{Preliminaries on Extremum Seeking Control}\label{sec:es:pre}


Extremum seeking (ES) control  is a type of model-free control  that uses only output feedback to
steer a dynamical system to a state where the output function attains an extremum \cite{ariyur2003real}. 
Therefore, ES can be interpreted as a continuous-time zeroth-order method to solve optimization problems, which essentially estimates  the gradient of the objective function based on exploratory probing signals.

\begin{figure}
    \centering
     \includegraphics[scale=0.358]{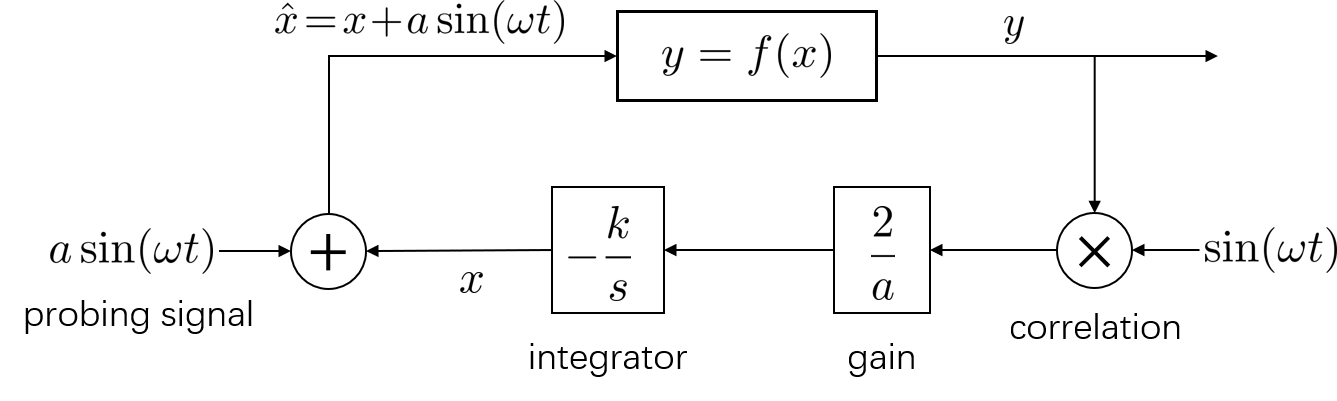}
    \caption{The block diagram of a simple ES scheme for solving $\min_x f(x)$.}
    \label{fig:es}
\end{figure}

Consider the problem of solving  $\min_x f(x)$. A straightforward idea is to employ the gradient descent dynamics, i.e., $\dot{x} = -k\cdot \nabla f(x)$.  However, this dynamical system is not implementable when the gradient $\nabla f$ or the mathematical form of $f$ is unknown.  To address this issue, ES control estimates the gradient $\nabla f(x)$ based on sinusoidal probing signals. 
The simplest ES scheme that  consists of  necessary components  is shown as 
Figure \ref{fig:es}. 
Starting from  state $x$, a sinusoidal probing signal $a\sin(\omega t)$ with  frequency $\omega$ and amplitude $a$ is added to $x$. Then the perturbed input $\hat{x}$ is fed into the static map $y=f(x)$, and the output $y$ is multiplied by the sinusoidal signal $\sin(\omega t)$, leading to 
$ f(x + a\sin(\omega t))\sin(\omega t)$. 
The control loop  is closed through the gain $\frac{2}{a}$ and the integrator $\frac{-k}{s}$. Thus the dynamics of this closed-loop feedback ES system can be formulated as 
\begin{align}\label{eq:simes}
    \dot{x} = -k \cdot \frac{2}{a}  f(x+a\sin(\omega t)) \sin(\omega t),
\end{align}
where $(a,\omega,k)$ are design parameters. 

We first state the fact  that the ES dynamics \eqref{eq:simes} 
with  small  $a$ and large $\omega$  behaves, approximately, like the gradient descent dynamics  $\dot{x}=-k\cdot\nabla f(x)$, which can steer  $x$ to a (local) minimum $x^*=\arg\min_{x}f(x)$ under appropriate conditions on $f(\cdot)$. We also note that to implement the ES dynamics \eqref{eq:simes}, one does not need the knowledge of function $f$ but only its measurement.

The rationale behind is that   for sufficiently large value of $\omega$, the ES dynamics \eqref{eq:simes}  exhibits a timescale separation property, where the fast time variation  is caused by the sinusoidal signal $\sin(\omega t)$, while the slow variation that is governed by the gain $k$  dominates the evolution of $x$. 
By averaging theory, one can obtain a time-invariant average dynamics that describes the main  trend of the evolution  of $x$. With small value of $a$, we  consider the following Taylor expansion  in the scalar case: $$ f(x+a\sin(\omega t)) =  f(x)
+ a\sin(\omega t)\frac{\partial f(x)}{\partial x} +   \mathcal{O}(a^2)$$ 
%
%
Thus the average dynamics of (\ref{eq:simes}) is given by 
\begin{align}\label{eq:avedyn}
 \dot{x}& =  -k\cdot h_{\mathrm{av}}(x) = -k\cdot \frac{\partial f(x)}{\partial x} +\mathcal{O}(a),
 \end{align}
 where 
 $$
 h_{\mathrm{av}}(x)\!:=\! 
 \frac{1}{T} \int_{0}^T \!\frac{2}{a}  f(x+a\sin(\omega t))\sin(\omega t)  dt
      =\! \frac{\partial f(x)}{\partial x}\! +\mathcal{O}(a)
$$
and $T=\frac{2\pi}{\omega}$. The average dynamics \eqref{eq:avedyn} is indeed the gradient descent flow plus a small perturbation  $\sO(a)$.
The same idea can be applied to the multivariate case with an appropriate choice of the (vector) frequencies $\omega$.

The above simple  case  explains the  basic principle of ES control. 
While 
a  practical ES problem can be much more complex, e.g., involving a plant dynamics,  multiple-input and multiple-output, high-pass/low-pass filters, etc. 
See \cite{ariyur2003real} for  a detailed introduction.
\begin{remark} \normalfont
The ES system \eqref{eq:simes} is somehow analogous to the  single-point zeroth-order iterative method \cite{flaxman2004online}, given by:
\begin{align}\label{eq:single}
    x_{k+1} = x_{k} - \eta \frac{1}{r} f(x+ru)u
\end{align}
where $k$ is the iteration number,
$\eta$ is the step size, $r$ is the  smoothing radius, and $u$ is a random sample from an exploratory distribution,  e.g.,  Gaussian, with zero mean. See \cite{Manzie09,Teel01,Khong13,PovedaTAC17} for
more studies on the connection between ES control and zeroth-order optimization methods.

\end{remark}
%




\section{Algorithm Design} \label{sec:algorithm}
In this paper, we aim to design a real-time voltage control algorithm that satisfies the following four  requirements:
\begin{itemize}
    \item [1)] \textbf{Asymptotic voltage limits}. Once a disturbance occurs, the controller can drive the monitored voltage magnitudes 
    $(v_j)_{j\in\sM}$ back to the acceptable interval  $[\underline{v}_j,\bar{v}_j]$.
    \item [2)]  \textbf{Hard capacity constraints}. 
    The power injection $\bx_i$ of the controllable device $i\in\sC$ should satisfy the physical power capacity constraints $\sX_i$ at all times.
    \item [3)] \textbf{Optimality}. The controllable devices are regulated in  an economically efficient way that minimizes the total operating cost.
    \item [4)] \textbf{Model-free}. Information of the power network (topology and line parameters), loads and other power injections is not required. 
\end{itemize}


In this section, we first solve the OVC model (\ref{eq:ovc}) with the projected primal-dual gradient dynamics, so that the solution dynamics can be interpreted as the voltage controller which meets the first three requirements above. Then we take the fourth requirement into account and develop a model-free  voltage control algorithm based on ES control.

\subsection{Projected Primal-Dual Gradient Dynamics}

We make the following two standard assumptions  on the OVC model (\ref{eq:ovc}) to render it a convex optimization problem with strong duality. We emphasize that these assumptions are mainly for theoretical analysis, and 
the proposed control algorithm  can be applied
to power systems with a
nonlinear power flow model, which is validated by our simulations.

\begin{assumption} \label{ass:con_sm}
 For all $i\in\sC$, the
 function $c_i(\cdot)$ is convex and has locally Lipschitz gradients, and the set $\mathcal{X}_i$ is  closed and convex.
 Also, the function 
 $v_j(\cdot)$ is affine for all $j\in\sM$.
\end{assumption}

%
\begin{assumption} \label{ass:finite}
The OVC problem (\ref{eq:ovc}) has a finite optimum, and the Slater's conditions hold for the problem (\ref{eq:ovc}). 
\end{assumption}
%

We  employ   the  projected primal-dual gradient dynamics (\textbf{P-PDGD}) method to solve the OVC model (\ref{eq:ovc}).
With  dual variables $\bla^+\!:=\!(\lambda_j^+)_{j\in\sM},\bla^-\!:=\!(\lambda_j^-)_{j\in\sM}$,
the saddle point problem of the OVC model (\ref{eq:ovc}) is  formulated as
\begin{align} \label{eq:saddle}
\begin{split}
    &\max_{\bla\geq 0}  \min_{\bx\in\mathcal{X}} \, L(\bx, \bla) =\sum_{i\in\mathcal{C}} c_i(\bx_i)\\
      &\quad +\sum_{j\in\sM}\Big[ \lambda_j^+(v_j(\bx) - \bar{v}_j) + \lambda_j^-(\underline{v}_j-v_j(\bx) )  \Big]
\end{split}
\end{align}
where $\bla:=[\bla^+;\bla^-]$ and $ L(\bx, \bla)$ denotes the Lagrangian function. Then we  solve problem (\ref{eq:saddle}) with  P-PDGD (\ref{eq:ppdgd}): 
\begin{subequations} \label{eq:ppdgd}
\begin{align}
   & \dot{\bx}_i = k_x\Big[ \mathrm{Proj}_{\sX_i}\big( \,\bx_i\, - \alpha_x\frac{\partial L(\bx,\bla)}{\partial \bx_i }    \big) -\,\bx_i \,\Big],\, i\in\sC \label{eq:ppdgd:x}\\
   & \dot{\lambda}_j^+ \! =k_\lambda\Big[ \mathrm{Proj}_{\R_+}\!\big( \lambda_j^+ + \alpha_\lambda \frac{\partial L(\bx,\bla)}{\partial \lambda_j^+ }
     \big) -\!\lambda_j^+ \Big],\, j\!\in\!\sM \\
& \dot{\lambda}_j^-  \!=k_\lambda \Big[ \mathrm{Proj}_{\R_+}\!\big( \lambda_j^- + \alpha_\lambda\frac{\partial L(\bx,\bla)}{\partial \lambda_j^- } \big) -\!\lambda_j^- \Big],\, j\!\in\!\sM
\end{align}
\end{subequations}
where  $k_x, k_\lambda, \alpha_x, \alpha_\lambda$ are positive parameters, and
the Lipschitz projection operator $\mathrm{Proj}_{\sX}(\cdot)$ is defined as
\begin{align}
    \mathrm{Proj}_{\sX}(\bx):= \underset{\by \in \sX}{\mathrm{argmin}}  \, ||\by -\bx||. 
\end{align}
The gradients in (\ref{eq:ppdgd}) are given by
\begin{subequations}
\begin{align}
   \frac{\partial L(\bx,\bla)}{\partial \bx_i }   & = \nabla c_i(\bx_i)
 + \sum_{j\in\sM}(\lambda_j^+-\lambda_j^-)\frac{\partial v_j(\bx)}{\partial \bx_i} \label{eq:grad:x}\\
\frac{\partial L(\bx,\bla)}{\partial \lambda_j^+ } & = v_j(\bx) -\bar{v}_j \label{eq:grad:la+}\\
\frac{\partial L(\bx,\bla)}{\partial \lambda_j^- }  & = \underline{v}_j-v_j(\bx). \label{eq:grad:la-}
\end{align}
\end{subequations}
$\blacktriangle$ Denote $\bz:= [\bx;\bla]$  and 
define $\mathcal{Z}:= \mathcal{X}\times \R_+^{2|\mathcal{M}|}$ as the feasible set  of $\bz$ in (\ref{eq:ppdgd}). 

\begin{remark} \label{remark:project}
(Projection of Dynamical System) \normalfont The
 projection method used in (\ref{eq:ppdgd}) is referred as  \emph{global projection} \cite{gao2003exponential}. By \cite[Lemma 3]{gao2003exponential},  it ensures that $\bz(t)\in \sZ$ for all time $t\geq 0$ when the initial condition $\bz(0)\in \sZ$.  For example, consider the dynamics of $\bx_i$. The intuition of this type of projection is that (\ref{eq:ppdgd:x}) attempts to take a step forward with stepsize $\alpha_x$ along the gradient descent direction, then checks whether the arrival point $\bx_i- \alpha_x\frac{\partial L(\bx,\bla)}{\partial \bx_i }$ is feasible to $\sX_i$. If feasible, (\ref{eq:ppdgd:x}) reduces to the ordinary gradient descent dynamics $\dot{\bx}_i \!=\! -k_x\alpha_x \frac{\partial L(\bx,\bla)}{\partial \bx_i }$, otherwise a projection is performed to guarantee the feasibility of $\bx_i$. Note that the P-PDGD \eqref{eq:ppdgd} is Lipschitz continuous; it differs from other types of discontinuous projections considered in literature, e.g.,  \cite{nagurney2012projected,8571158,9147814}, which project the dynamics onto the tangent cone of the feasible set, and thus they need the sophisticated analysis tools  for discontinuous  
dynamical systems.
\end{remark}


%

Following the P-PDGD (\ref{eq:ppdgd}), the state $\bx(t)$ will remain within the feasible set $\sX$ and
converge  to  a  steady-state operating point that is an optimal solution of the OVC problem (\ref{eq:ovc}). 
This is restated formally as  Theorem \ref{thm:ppdgdcon}. 
\begin{theorem} \label{thm:ppdgdcon}
(Global Asymptotical Stability.) 
Under  Assumption \ref{ass:con_sm} and \ref{ass:finite}, 
with initial condition $\bz(0)\in \sZ$, the trajectory $\bz(t)$ of the P-PDGD (\ref{eq:ppdgd}) will stay within $\sZ$ for all $t\geq 0$ and globally asymptotically converge to an optimal solution $\bz^*:=[\bx^*;\bla^*]$ of the saddle point problem (\ref{eq:saddle}), 
 %
where $\bx^*$ is an optimal solution of the OVC problem \eqref{eq:ovc}.
\end{theorem}
The proof of Theorem  \ref{thm:ppdgdcon} mainly follows the asymptotical stability of globally projected (primal-dual) dynamical systems \cite[Lemma 2.4]{bansode2019exponential} \cite{gao2003exponential}. A detailed proof is provided in Appendix \ref{sec:app1}.
As a result, the P-PDGD (\ref{eq:ppdgd}) can be regarded as the voltage control mechanism that meets the first three requirements above. In the next subsection, we will take into account the fourth requirement and develop a model-free   control algorithm based on the proposed P-PDGD (\ref{eq:ppdgd}).
\subsection{Model-Free Voltage Control Algorithm}

The P-PDGD (\ref{eq:ppdgd}) cannot be  implemented without knowledge of the system model  $\bv(\bx)$. Note that there are two occasions in the P-PDGD (\ref{eq:ppdgd}) where this model is needed: 
\begin{itemize}
    \item [1)] The gradients $\frac{\partial v_j(\bx)}{\partial \bx_i}$ in (\ref{eq:grad:x}) for $i\in\sC, j\in\sM$;
       \item [2)] The functions $v_j(\bx)$ in (\ref{eq:grad:la+}) (\ref{eq:grad:la-}) for $j\in\sM$.
\end{itemize}
To develop a model-free controller, we propose the following two strategies accordingly:

\vspace{0.1cm}
\textbf{Strategy 1)}:  Use ES control  to ``estimate" the gradients $\frac{\partial v_j(\bx)}{\partial \bx_i}$ 
for all $i\in\sC, j\in\sM$.

\vspace{0.1cm}
\textbf{Strategy 2)}: Substitute the function $v_j(\cdot)$ by the real-time voltage measurement $v_j^{\mathrm{mea}}(t)$ for all $j\in\sM$. 

\vspace{0.1cm}
To implement Strategy 1), we add a small sinusoidal probing signal to each power injection with
\begin{align} \label{eq:xhat}
    \hat{\bx}_i(t) = \bx_i(t) + a \sin(\bom_i t),\quad i\in\sC
\end{align}
where $a$ is the small amplitude\footnote{For notational simplicity, we adopt an identical amplitude $a$ for all power injections here. In practice, different amplitude parameters can be used.} and the sinusoidal signal is 
\begin{align}
    \sin(\bom_i t):= \begin{cases}
   \ \sin(\omega_i t), & i\in\sC_{\mathrm{svc}}\\
    [ \sin(\omega_i^p t), \sin(\omega_i^q t)]^\top, &i\in \sC_{\mathrm{dg}}.
    \end{cases}
\end{align}
$\blacktriangle$ Let $N= |\sC_{\mathrm{svc}}|+2|\sC_{\mathrm{dg}}|$ be the dimensionality of the decision variable $\bx$. 
Define $\sin(\bom t):=(\sin(\bom_i t))_{i\in\sC}\in\R^{N}$ as the column vector that collects all the  sinusoidal signals. The frequencies $\bom:=(\bom_i)_{i\in\sC}$ are selected as 
\begin{align} \label{eq:fredef}
  \qquad\qquad  \omega_n =\frac{2\pi}{\varepsilon_\omega}\kappa_n,\quad \forall n\in [N]:=\{1,\cdots,N\}
\end{align}
where $\varepsilon_\omega$ is a small positive parameter and $\kappa_i\neq \kappa_j$ for all $i\neq j$ are rational numbers. In this way, each element $x_n$ in $\bx$ is assigned with a particular frequency $\omega_n$.

Based on the above description, the P-PDGD (\ref{eq:ppdgd}) is modified as the  \textbf{ES-P-PDGD} (\ref{eq:es:ppdgd}):
\begin{subequations} \label{eq:es:ppdgd}
\begin{align}
    \dot{\bx}_i \!&= \!k_x\Big[ \mathrm{Proj}_{\hat{\sX}_i}\big( \bx_i \!-\! \alpha_x \bh_i(\bx_i,\bla,\bxi_i^j)  \big)\! -\!\bx_i \Big],\ i\in\sC \label{eq:es:x}\\
    \dot{\lambda}_j^+\! & \!=k_\lambda \Big[ \mathrm{Proj}_{\R_+}\!\big( \lambda_j^+ + \alpha_\lambda (\mu_j -\bar{v}_j)
     \big) \!-\!\lambda_j^+ \Big], \ \, j\!\in\!\sM\label{eq:es:la+} \\
 \dot{\lambda}_j^-\! &\! =k_\lambda \Big[ \mathrm{Proj}_{\R_+}\!\big( \lambda_j^- + \alpha_\lambda(\underline{v}_j - \mu_j) \big)\! -\!\lambda_j^- \Big],\ \, j\!\in\!\sM \label{eq:es:la-}\\
      \dot{\bxi}_i^j &= \frac{1}{\epsilon} \Big[ \! -{\bxi}_i^j + \frac{2}{a} {v}_j (\hat{\bx}(t))\sin(\bom_i t) \Big],\ j\!\in\!\sM, i\in\sC\label{eq:es:xi}
\\
 \dot{\mu}_j &= \frac{1}{\epsilon}\Big[ -\!\mu_j + v_j (\hat{\bx}(t)) \Big], \qquad\qquad\qquad\quad\   j\!\in\!\sM \label{eq:es:mu}
\end{align}
\end{subequations}
where $\epsilon$ is a small positive parameter, and
\begin{subequations}
\begin{align}
\hat{\bx}(t)& := \bx(t) + a\sin(\bom t)\\
   \bh_i(\bx_i,\bla,\bxi_i^j) &:= \nabla c_i(\bx_i) + \sum_{j\in\sM}(\lambda_j^+-\lambda_j^-)\bxi_i^j.
\end{align} 
\end{subequations}
The key difference between P-PDGD (\ref{eq:ppdgd}) and ES-P-PDGD (\ref{eq:es:ppdgd}) is the introduction of new variables $\bxi\!:=\!(\bxi_i^j)_{j\in\sM,i\in\sC}$ and $\bmu\!:=\!(\mu_j)_{j\in\sM}$. We explain the rationale and benefits of this modification with  Remark \ref{remark:muxi}. To ensure the actual power injection $\hat{\bx}_i\in \sX_i$, 
we replace $\sX_i$ with the shrunken feasible set $\hat{\sX}_i$ \eqref{eq:hatX} in  (\ref{eq:es:x}). As $a\to 0^+$, 
 $\hat{\sX}$   recovers to  $\sX$.
\begin{align} \label{eq:hatX}
  \hat{\mathcal{X}}_i\!:= \!
  \begin{cases}
 \underline{q}_i\!+\!a\!\leq q_i\leq \bar{q}_i\!-\!a, \qquad i\in\mathcal{C}_{\mathrm{svc}}  \\
   \underline{p}_i\!+\!a\!\leq p_i\leq \bar{p}_i\!-\!a,\,
   p_i^2\!+\!q_i^2\!\leq\! (\bar{s}_i\!-\!\sqrt{2}a)^2 , \, i\!\in \mathcal{C}_{\mathrm{dg}}. 
  \end{cases}
\end{align}

\begin{remark}\label{remark:muxi}
(Fast Dynamics of $\bxi$ and $\bmu$.) \normalfont In essence, $\bxi_i^j$ and $\mu_j$ are the real-time approximations of the gradient $\frac{\partial {v}_j}{\partial \bx_i}$ and the value $v_j$, respectively. 
The intuition behind is that by setting $\epsilon$ sufficiently small, the dynamics of $\bxi^j_i$ and $\mu_j$, i.e., (\ref{eq:es:xi}) (\ref{eq:es:mu}), operate in a faster time scale compared to the dynamics of $(\bx,\bla)$.  The advantages of introducing these fast dynamics of $\bxi$ and $\bmu$  include:
\begin{itemize}
    \item [1)]  It facilitates the analysis of the algorithm via averaging theory, since the time-varying sinusoidal signals do not appear inside the projection operators.
    Moreover, the fast dynamics are linear, which can be easily handled by singular perturbation theory. 
    \item [2)] The fast dynamics (\ref{eq:es:xi}) (\ref{eq:es:mu}) can be seen as low-pass filters, which can diminish the oscillations and improve the transient performance of the closed-loop system.
\end{itemize}
\end{remark}
Since $\hat{\bx}(t)$  is the actual power injection to the physical system at time $t$, we can substitute $v_j(\hat{\bx})$ with the voltage  measurement $v_j^{\mathrm{mea}}(t)$ in (\ref{eq:es:ppdgd}), i.e., Strategy 2). 
 Consequently,
 we develop the   model-free optimal voltage control (\textbf{MF-OVC}) algorithm as Algorithm \ref{alg:voltctrl}, which is indeed the ES-P-PDGD \eqref{eq:es:ppdgd} with the measurement substitution.
 \begin{algorithm}
 \caption{Model-Free Optimal Voltage Control (MF-OVC) Algorithm.}
 \begin{algorithmic}  \label{alg:voltctrl}
\STATE At every time $t$, perform the following steps:
 \STATE $\bullet$ \textbf{Each monitored bus} $j\in\sM$ measures the local voltage magnitude $v^{\mathrm{mea}}_j(t)$, updates $( \lambda_j^+,\lambda_j^-,\mu_j)$ according to 
  \begin{subequations}
   \begin{align}
      & \text{Equations (\ref{eq:es:la+}) (\ref{eq:es:la-})} \\
       & \dot{\mu}_j = \frac{1}{\epsilon}\Big[ -\mu_j + v^{\mathrm{mea}}_j(t) \Big]
  \end{align}
  \end{subequations}
 and broadcasts $(v^{\mathrm{mea}}_j(t), \lambda_j^+(t),\lambda_j^-(t))$ to every controllable device $i\in \sC$.
 \STATE   $\bullet$  \textbf{Each controllable device} $i\in\sC$ updates $(\bx_i,\bxi_i^j)$ by 
 \begin{subequations}
  \begin{align}
   & \text{Equation (\ref{eq:es:x})}\\
       & \dot{\bxi}_i^j  = \frac{1}{\epsilon} \Big[  \!-{\bxi}_i^j  + \frac{2}{a}{v}_j^{\mathrm{mea}}(t)\sin(\bom_i t)  \Big],\ j\!\in\!\sM
  \end{align}
 \end{subequations}
  and executes  power injection $\hat{\bx}_i(t)= \bx_i(t)+ a\sin(\bom_i t)$.

 \end{algorithmic} 
 \end{algorithm}

%
\begin{figure}
    \centering
    \includegraphics[scale=0.38]{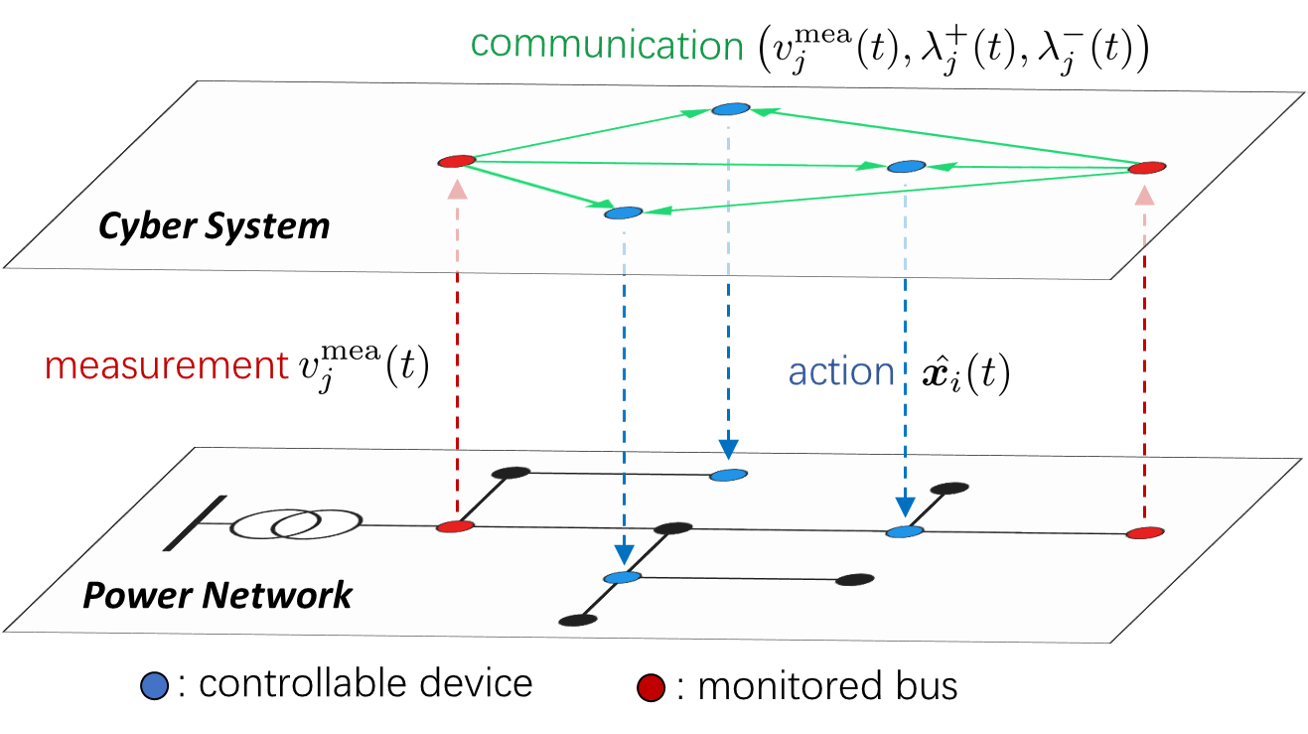}
    \caption{Schematic of  the proposed MF-OVC mechanism.}
    \label{fig:voltcon}
\end{figure}
The implementation of the proposed MF-OVC algorithm is illustrated in Figure \ref{fig:voltcon}.  
 Each monitored bus $j\in\sM$ measures its local voltage magnitude $v_j^{\mathrm{mea}}$  from the physical  layer, then updates $(\mu_j, \lambda_j^+,\lambda_j^-)$ and communicates $(v^{\mathrm{mea}}_j, \lambda_j^+,\lambda_j^-)$ in the cyber  layer. Each controllable device $i\in\sC$ updates $(\bxi_i^j,\bx_i)$ based on the received information, and the power injection command $\hat{\bx}_i$ is  executed  in the physical layer.  Then the power network responses to the power injection $\hat{\bx}$ and presents the corresponding voltage profiles $\bv(\hat{\bx})$. This forms a closed-loop feedback control system. Although the MF-OVC  algorithm is  developed based on a static OVC problem (\ref{eq:ovc}), it can   adapt fast to dynamical      system environments and      handle voltage violation under  time-varying  power disturbances, due to the feedback mechanism and exploitation of  real-time measurements.
   This is validated by the simulations in Section \ref{sec:sim:timevary}.
As a result, the proposed  algorithm unifies all the merits described in the introduction section.

\section{Performance Analysis} \label{sec:performance}
This section presents the theoretical analysis on the performance of the proposed MF-OVC algorithm. In particular, we focus on the ES-P-PDGD  (\ref{eq:es:ppdgd}) and study its stability  properties as well as its robustness to measurement noises.

\subsection{Stability Analysis of ES-P-PDGD }

\noindent $\blacktriangle$ Denote $\bz\!:= \![\bx;\bla]$ 
and 
$\hat{\sX}:= \prod_{i\in\sC} \hat{\sX}_i$. Let 
$\hat{Z}:= \hat{\sX}\times \R^{2|\mathcal{M}|}_+$ be the feasible set of $\bz$ in the ES-P-PDGD (\ref{eq:es:ppdgd}),  and $K := (2|\sC_{\mathrm{dg}}| \!+\!|\sC_{\mathrm{svc}}|\!+\!1)|\sM|$ be the dimensionality of $[\bxi;\bmu]$.  
Denote $\hat{\sA}$ as the saddle point set  for the saddle point problem \eqref{eq:saddle} with  $\hat{\sX}$, i.e., any point $\hat{\bz}^*\in \hat{\sA}$ is an optimal solution of \eqref{eq:saddle} with  $\hat{\sX}$.\footnote{Here, the notations with  ``$\wedge$" on the head represent the counterparts with the feasible set $\hat{\sX}$.} Denote the distance between $\bz$ and $\hat{\sA}$ as 
\begin{align*}
    ||\bz||_{\hat{\sA}}:= \underset{\bm{\alpha}\in \hat{\sA}}{\inf}\, ||\bz-\bm{\alpha}||.
\end{align*}

\begin{definition} \normalfont
A continuous function $\beta(r,t): \R_+\times \R_+\to \R_+$ is said to be of class-$\mathcal{KL}$ if it is zero at zero and strictly increasing in the first argument $r$, and non-increasing in the second argument $t$ and converging to zero  as  $t\to +\infty$.
\end{definition}

The  stability of ES-P-PDGD (\ref{eq:es:ppdgd}) is stated as Theorem \ref{thm:spas}. 



\begin{theorem}\label{thm:spas}
(Semi-Global Practical Asymptotical Stability.) Suppose that the saddle point set $\hat{\sA}$ is compact.
Under  Assumption \ref{ass:con_sm} and \ref{ass:finite}, 
 there exists a class-$\mathcal{KL}$ function $\beta$ such that for any compact set $\mathcal{D}\subset \hat{Z}\times \R^{K}$ of initial condition, 
and any desired precision $\nu>0$, there exists $\epsilon^*>0$ such that for any $\epsilon\in(0,\epsilon^*)$, there exists $a^*>0$ such that for any $a\in(0,a^*)$, there exists $\varepsilon_\omega^*>0$ such that for any $\varepsilon_\omega\in(0,\varepsilon_\omega^*)$, the trajectory $\bz(t)$ of the  ES-P-PDGD (\ref{eq:es:ppdgd}) satisfies 
\begin{align}\label{eq:spas}
    ||\bz(t)||_{\hat{\sA}}\leq \beta(||\bz(0)||_{\hat{\sA}},\,  t) +\nu, \ \ \forall t\geq 0.
\end{align}
\end{theorem}

We prove Theorem \ref{thm:spas} using averaging theory and singular perturbation theory \cite{khalil2002nonlinear,Wang:12_Automatica}.  The  detailed proof of Theorem \ref{thm:spas} is provided in Appendix \ref{sec:app2}. 



\begin{remark} \normalfont
We explain the key observations of Theorem \ref{thm:spas} as follows:

$\bullet$ Due to the small probing sinusoidal signals $a\sin(\bom t)$ in the ES-P-PDGD (\ref{eq:es:ppdgd}), the state $\bz$ will not converge to a fixed point anymore,
but rather to a small $\nu$-neighborhood of $\hat{\sA}$.
This property is described by the bound (\ref{eq:spas}).
By setting the parameters $(\epsilon, a, \varepsilon_\omega)$ sufficiently small, one can make this precision $\nu$
 as small as desired.


$\bullet$ As $(\epsilon,a,\varepsilon_{\omega})\to0^+$, the
ES-P-PDGD \eqref{eq:es:ppdgd} recovers the same convergence rate of the P-PDGD \eqref{eq:ppdgd},  as indicated in the proof of Theorem \ref{thm:spas}.

$\bullet$ As stated in Theorem \ref{thm:spas}, the tuning order of parameters  is relevant: first set $\epsilon$ sufficiently small, then $a$, and lastly $\varepsilon_\omega$. This order comes mainly  from the proof and can guide us on how to tune these parameters in practice.



%
%

\end{remark}

\begin{remark}
\normalfont We note that the assumption of a compact saddle point set $\hat{\sA}$ in Theorem \ref{thm:spas} is standard for the application of averaging theory and singular perturbation theory.  For the OVC problem \eqref{eq:ovc}, if the cost function $c_i(\cdot)$ is strictly convex for all $i\in\sC$ and the Jacobian matrix $\nabla_{\bx} \bv(\bx)$ is of full row rank, one can prove that the saddle point set $\hat{\sA}$ is singleton, i.e., the optimal solution  of the saddle point problem \eqref{eq:saddle} is \emph{unique}, by \cite[Proposition 1]{qu2018exponential}. In practice, the condition that the Jacobian matrix $\nabla_{\bx} \bv(\bx)$ is of full row rank can be satisfied when the  number of controllable devices is more than the monitored buses in the distribution system \cite{6760555}. 
\end{remark}


\subsection{Robustness to Measurement Noise}

The proposed  algorithm purely relies on the voltage measurement for control. Accordingly, the following
corollary of Theorem \ref{thm:spas} \cite{PovedaNaLi2019} indicates that this algorithm is robust to \emph{small} additive state measurement noise. Moreover, the numerical simulations in Section \ref{sec:noise}  verify the robustness even when the noise is relatively large. 
\begin{corollary} \label{thm:robust}
(Structural Robustness.) For any tuple of $(\epsilon,a,\varepsilon_\omega)$ that induces the bound \eqref{eq:spas}, under the same conditions in Theorem \ref{thm:spas},
there exists $\rho^*>0$ such that for any  measurement noise $\bm{d}:\R_+\!\to\! \R^{|\sM|}$ with $\sup_{t\geq 0} ||\bm{d}(t)||\leq \rho^*$, the trajectory $\bz(t)$ of the ES-P-PDGD (\ref{eq:es:ppdgd}) with additive state measurement noise $\bm{d}$ satisfies
\begin{align}\label{eq:robust}
    ||\bz(t)||_{\hat{\sA}}\leq \beta(||\bz(0)||_{\hat{\sA}},\,  t) +2\nu, \ \ \forall t\geq 0.
\end{align}
\end{corollary}
%
Comparing with (\ref{eq:spas}), the ES-P-PDGD \eqref{eq:es:ppdgd} with small additive measurement noise $\bd$  maintains similar convergence results,
and noise $\bm{d}$  leads to an additional $\nu$ term in \eqref{eq:robust}.
Besides, this robustness property can be extended to other small additive perturbations.

\section{Numerical Simulations} \label{sec:simulation}

In this section, we demonstrate the  performance of the proposed MF-OVC algorithm via numerical simulations. Specifically, we test the MF-OVC algorithm under  step and continuous power disturbances. The impact of noises in voltage measurements is studied  numerically as well.

\subsection{Simulation Setup}

\begin{figure}
    \centering
        \includegraphics[scale=0.525]{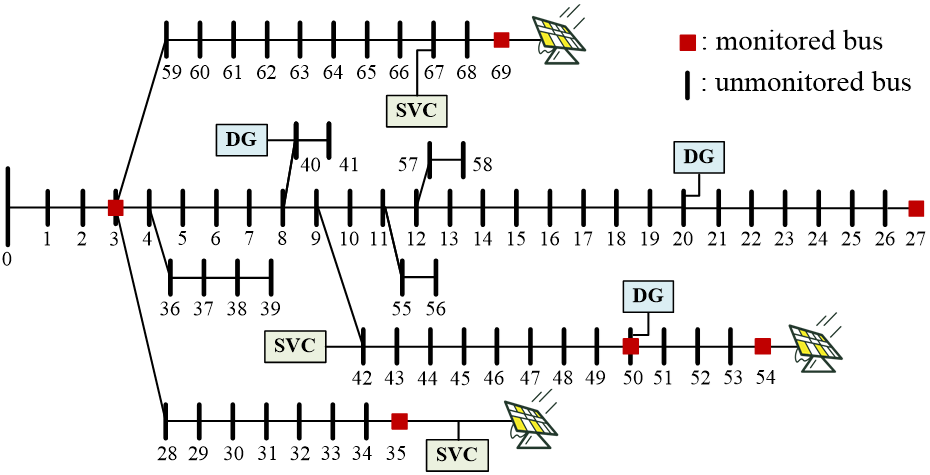}
    \caption{The modified PG\&E 69-bus distribution feeder.}
    \label{fig:PGE49}
\end{figure}

The modified PG\&E 69-bus distribution system, shown as Figure \ref{fig:PGE49}, is used as the test system. There are three PV plants at bus 35, 54 and 69, which operate in the maximum power point tracking mode. The controllable devices include three SVCs (located at bus 35, 42 and 67) and three DGs (located at bus 20, 40, 50). 
 Their power capacity limits  are set to
 \begin{align*}
    & \underline{q}_i = -1.5 \text{ MVar}, \, \bar{q}_i = 0.6 \text{ MVar}, \quad i\in\sC_{\mathrm{svc}} \\
     &\underline{p}_i = 0\text{ MW},\, \bar{p}_i = 1.5\text{ MW},\, \bar{s}_i = 1.8 \text{ MVA}, \quad i\in\sC_{\mathrm{dg}}.
 \end{align*}
We select bus 3, 27, 35, 50, 54 and 69 as the monitored buses. The voltage of bus 0 (slack bus) is 10.5 kV (1 p.u.), and the lower and upper bounds of voltage magnitude are set as 0.95 p.u. and 1.05 p.u., respectively. We use the quadratic cost function  (\ref{eq:cost}) with the coefficient $c_i^{\mathrm{svc}}=0.1, c_{p,i}^{\mathrm{dg}} =1, c_{q,i}^{\mathrm{dg}} = 0.5$.
For the MF-OVC algorithm, we set $a = 0.05$, $\epsilon = 0.02$,  $\varepsilon_\omega = 0.05$, and $\kappa_n = 2n-1$ for $n=1,\cdots,9$.

Although  an affine voltage function $\bv(\bx)$ is assumed for theoretical analysis, we perform all the simulations based on a \emph{full nonconvex AC power flow model} using the Matpower software \cite{matpower}. 
\subsection{Static Voltage Control Under Step Power Change} \label{sec:sim:step}
\begin{figure}
    \centering
     \includegraphics[scale=0.29]{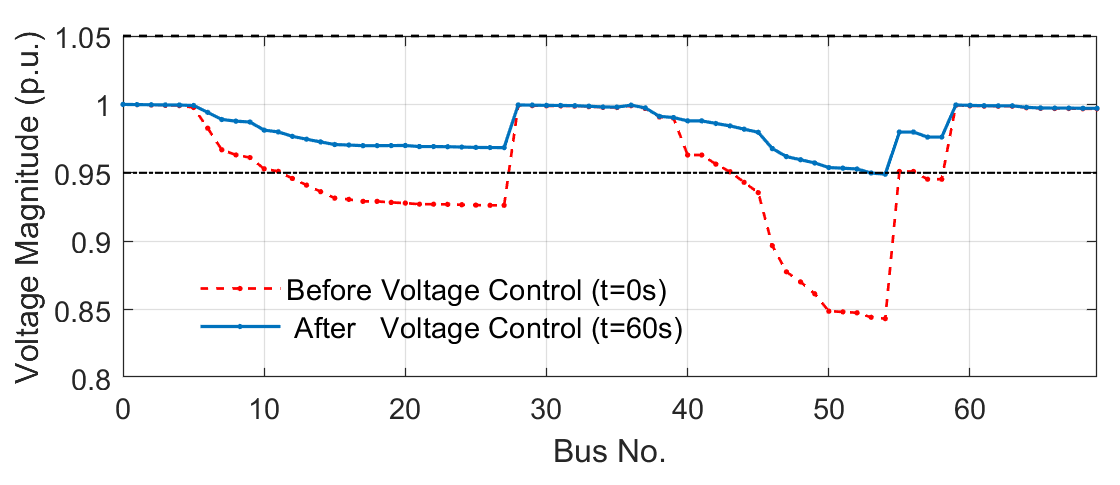}
    \caption{The voltage magnitude profiles before and after voltage control (black dotted lines: upper (1.05 p.u.) and lower (0.95 p.u.) voltage limits).}
    \label{fig:lowvolt}
\end{figure}

\begin{figure}
    \centering
    \includegraphics[scale=0.285]{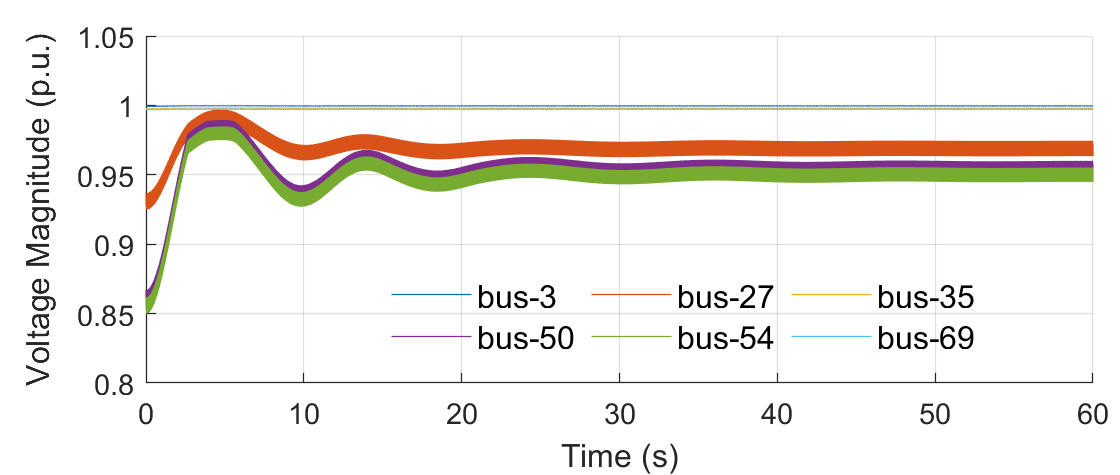}
    \caption{Voltage dynamics of the monitored buses under step power change.}
    \label{fig:stepvolt}
\end{figure}

\begin{figure}
    \centering
     \includegraphics[scale=0.29]{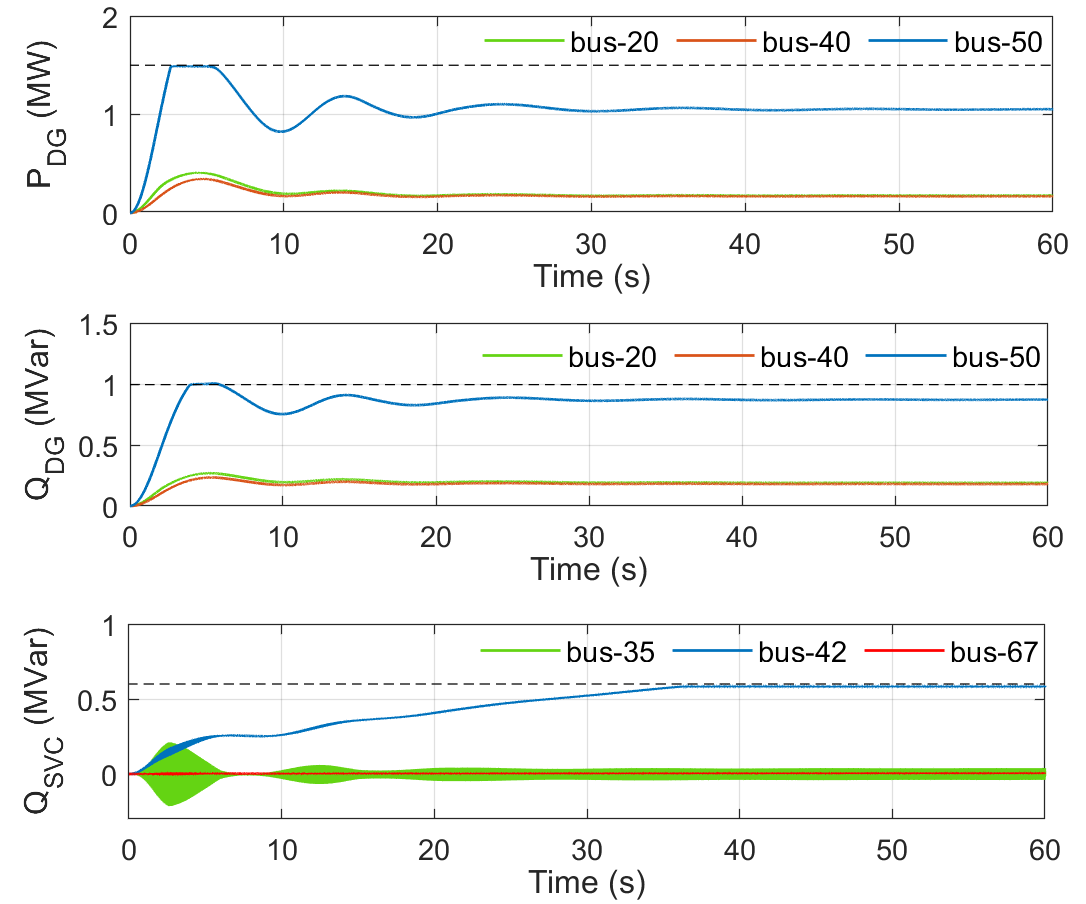}
    \caption{The active/reactive power outputs of DGs and SVCs (black dashed lines: the corresponding power capacity limits).}
    \label{fig:power}
\end{figure}

Consider the test scenario when the three PV plants are suddenly shut down at time $t=0$ and all loads remain fixed. Due to the curtailment of PV generation and heavy loads,   the 
 voltage profiles of the test system decrease to a low level. It  leads to  voltage violation at many buses,  shown as the red dashed curve in Figure \ref{fig:lowvolt}. 
 We ran the proposed MF-OVC algorithm for voltage regulation  from the start time $t=0$. 
 The simulation results are shown as Figure \ref{fig:stepvolt} and Figure \ref{fig:power}.

From Figure \ref{fig:stepvolt}, it is observed that the proposed MF-OVC algorithm can quickly bring  the voltage magnitudes of monitored buses back to the acceptable range. The small high-frequency oscillations in voltage are caused by the exploratory sinusoidal signals in the MF-OVC algorithm. As a result, the  voltage profiles of the entire test system were restored to the acceptable level (see the blue curve in Figure \ref{fig:lowvolt}),
due to the selection of representative monitored buses.
Figure \ref{fig:power} illustrates the dynamics of the power outputs of DGs and SVCs. It is seen that the power outputs converge to  fixed values (with small oscillations) within  tens of seconds, and the power capacity constraints are satisfied all the time. Besides, we  solve the OVC model (\ref{eq:ovc}) to obtain the optimal solution $\bx^*$\footnote{We solve the OVC model (\ref{eq:ovc}) with the CVX package \cite{cvx}, and the linearized Distflow model \cite{8779692} is used as the power flow model for $\bv(\bx)$.}, which turns out to be the converged  values in Figure \ref{fig:power}. It
 verifies the optimality of the MF-OVC algorithm. 


\subsection{Dynamic Voltage Control Under Continuous Change}
\label{sec:sim:timevary}

\begin{figure}
    \centering
    \includegraphics[scale=0.3]{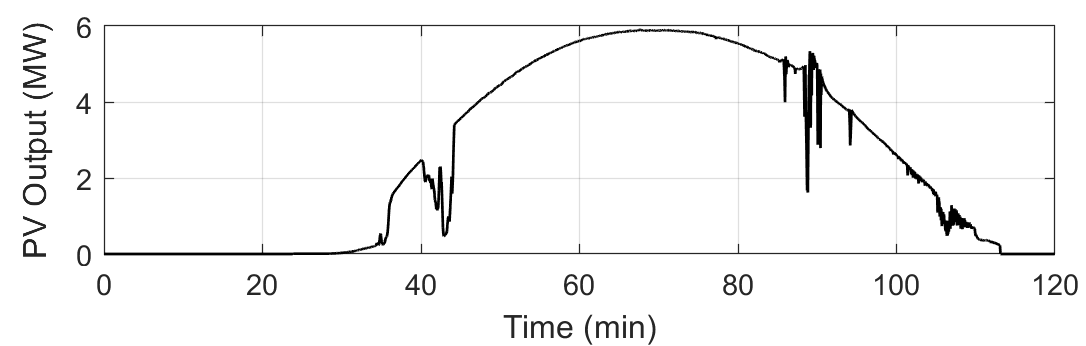}
    \caption{The time-varying total PV generation over two hours.}
    \label{fig:pvpower}
\end{figure}

We then test the performance of the proposed MF-OVC algorithm under time-varying loads and PV generations. We add a $10\%$ random perturbation to the total load,  and a real-world PV generation profile, shown as Figure \ref{fig:pvpower}, is applied to the three PV plants in the test system.  We ran the proposed MF-OVC algorithm for voltage regulation and compared it with the case without voltage control.  

\begin{figure}
    \centering
    \includegraphics[scale = 0.302]{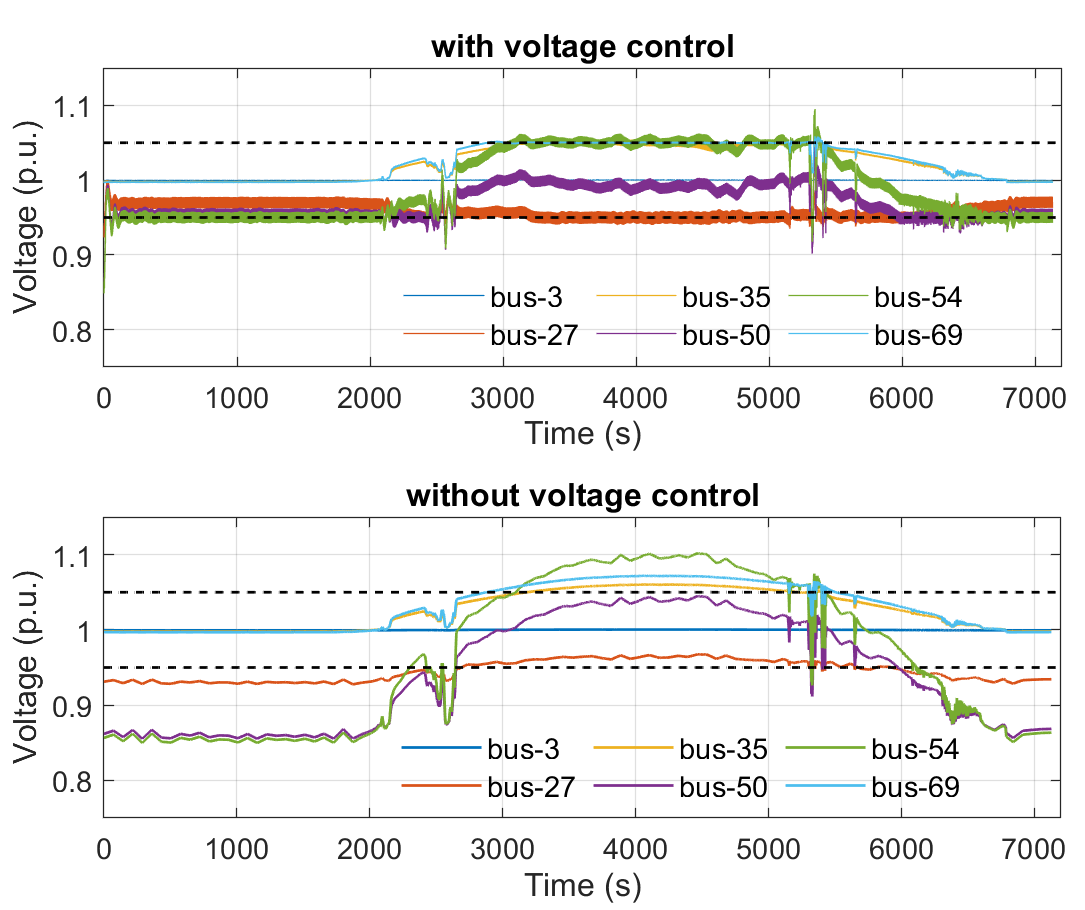}
    \caption{Voltage dynamics of the monitored buses under continuous power disturbances (black dashed lines: upper (1.05 p.u.) and lower (0.95 p.u.) voltage limits).}
    \label{fig:volcompare}
\end{figure}

The simulation results are illustrated in Figure \ref{fig:volcompare}. In the absence of voltage control, the test system violates the lower voltage limit (0.95 p.u.) when the PV generation is low, and the upper voltage limit (1.05 p.u.) when the PV generation is high. In contrast, the proposed MF-OVC algorithm can effectively adapt  to the  continuous power disturbances and maintain the voltage profiles within the acceptable range. 

\subsection{Impact of Measurement Noise}\label{sec:noise}
\begin{figure}
    \centering
    \includegraphics[scale=0.3]{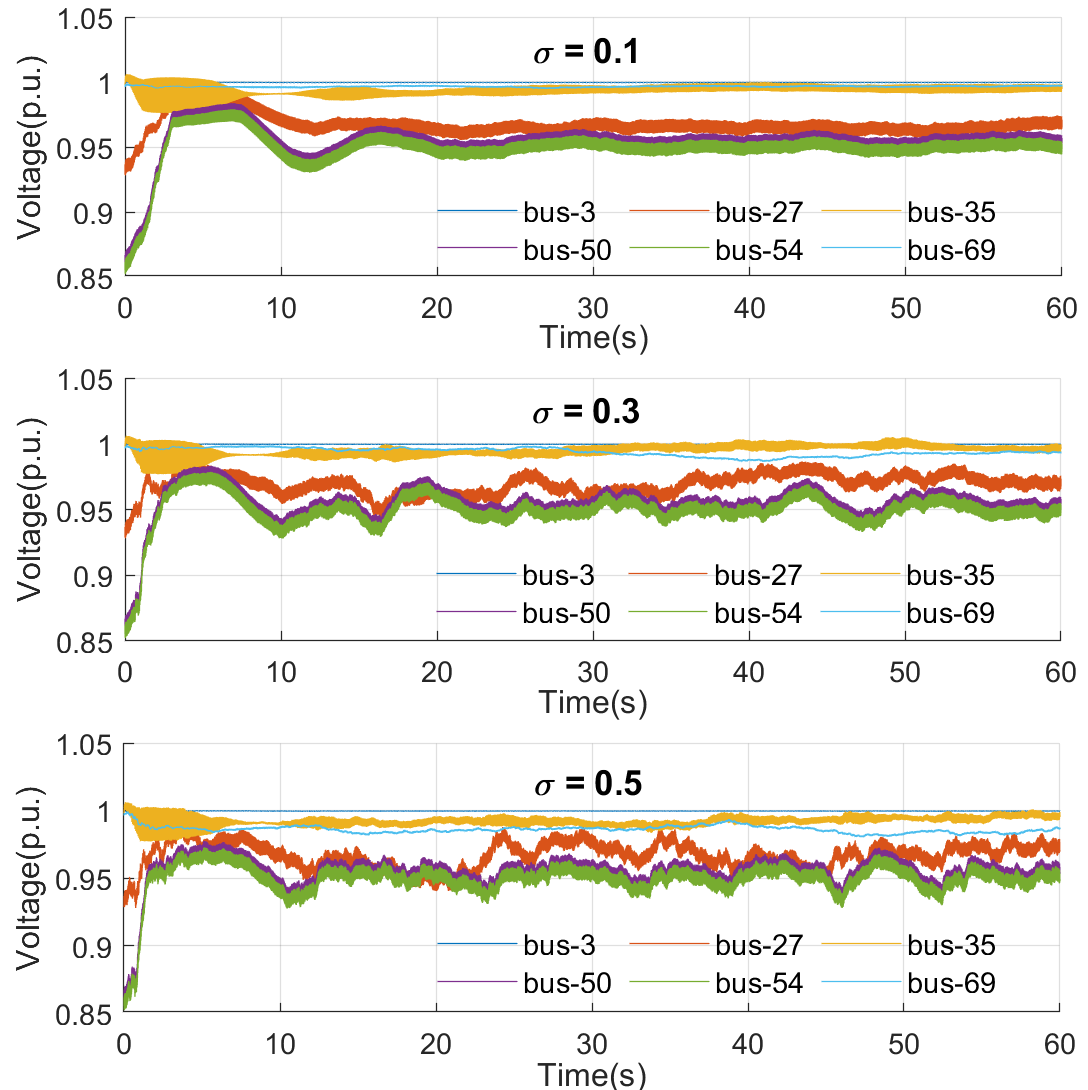}
    \caption{Voltage dynamics of the monitored buses with noisy voltage measurements.}
    \label{fig:noise}
\end{figure}

To study the impact of measurement noises, we consider the noisy voltage measurement $\tilde{v}^{\mathrm{mea}}_j(t)$, whose deviation from the base voltage value (1 p.u.) follows (\ref{eq:noise}):
\begin{align} \label{eq:noise}
   \tilde{v}^{\mathrm{mea}}_j(t) - 1 =  (v_j(\bx(t)) - 1)\times(1+\delta_j(t) )
\end{align}
where $v_j(\bx(t))$ denotes the true voltage magnitude, and 
 $\delta_j$ is the perturbation ratio. We assume that $\delta_j$ is a Gaussian random variable with $\delta_j\sim \mathcal{N}(0,\sigma^2)$, which is  independent across time $t$ and other monitored buses.  We tune the standard deviation $\sigma$ from 0.1 to 0.5 to simulate different levels of noises and test the performance of the MF-OVC algorithm under step power changes. The simulation results are shown as Figure \ref{fig:noise}, and the noiseless case with $\sigma =0$ is illustrated in Figure \ref{fig:stepvolt}. As expected, larger noise amplitudes lead to higher oscillations in the voltage dynamics. While
 the MF-OVC algorithm is robust to the voltage measurement noises and  can bring the voltage profiles back to the acceptable interval in all the cases.

\section{Conclusion}\label{sec:conclusion}

In this paper, we developed a real-time model-free  optimal voltage control algorithm based on projected primal-dual gradient dynamics and extremum seeking control. The proposed algorithm operates purely based on the voltage measurement and does not require any other network information. 
With appropriate parameters, this algorithm can effectively bring the monitored voltage magnitudes back to the acceptable range with minimum operational cost, while respecting the power capacity constraints all the time.
Numerical simulations on a modified PG\&E 69-bus distribution feeder demonstrated  that the proposed algorithm is capable of handling voltage violation under step or continuous power disturbances, and is robust to measurement noises.

\appendix
\section{}

\subsection{Proof of Theorem \ref{thm:ppdgdcon}} \label{sec:app1}
We first have the result of the following Proposition \ref{prop:equi} that connects the saddle point problem (\ref{eq:saddle}) and the P-PDGD (\ref{eq:ppdgd}). This proposition can be proved by  checking the KKT conditions of \eqref{eq:saddle} and using \cite[Theorem 3.25]{ruszczynski2011nonlinear}.
\begin{proposition}\label{prop:equi}
The optimal solutions of the saddle point problem (\ref{eq:saddle}) are equivalent to the equilibrium points of the P-PDGD (\ref{eq:ppdgd}). 
\end{proposition}

We then study the stability properties of the P-PDGD \eqref{eq:ppdgd}.
Denote $\bz:=[\bx;\bla]$, and let $\bz^*:=[\bx^*;\bla^*]$ be the optimal solution of \eqref{eq:saddle}. Define the mapping $$H(\bz):=  [ \nabla_{\bx} L(\bx,\bla); -\nabla_{\bla} L(\bx,\bla)],$$ where $L(\bx,\bla)$ is the Lagrangian function in \eqref{eq:saddle}. The P-PDGD \eqref{eq:ppdgd} can be written in  compact form\footnote{Without loss of generality, we set the time constant $k=1$ for simplicity.} as:
\begin{align} \label{eq:comp}
    \dot{\bz} = \mathrm{Proj}_{\sZ}(\bz - \alpha H(\bz)) -\bz:= \bm{f}(\bz),
\end{align}
where $\mathcal{Z}:= \mathcal{X}\times \R_+^{2|\mathcal{M}|}$. Since $\mathrm{Proj}_{\sZ}(\cdot)$ is a singleton and Lipschitz  on $\R^{|\sZ|}$ with constant $L = 1$ \cite[Proposition 2.4.1]{clarke1990optimization}, the dynamics $\bm{f}(\bz)$ in \eqref{eq:comp} is locally Lipschitz  on $\sZ$ by Assumption \ref{ass:con_sm}. Moreover, by \cite[Lemma 3]{gao2003exponential}, we have that $\bz(t)\in \sZ$ for all time $t\geq 0$ whenever $\bz(0)\in \sZ$. 

Next, consider the following Lyapunov function  $V$:
\begin{align} \label{eq:lyap}
    V(\bz) & := \frac{1}{2}||\bz-\bz^*||^2 + L(\bx,\bla^*) -L(\bx^*,\bla) \\
    & \ =  \frac{1}{2}||\bz-\bz^*||^2 + L(\bx,\bla^*)-L(\bx^*,\bla^*)\nonumber \\
    &\quad + L(\bx^*,\bla^*)-L(\bx^*,\bla) \nonumber \geq  \frac{1}{2}||\bz-\bz^*||^2.
\end{align}
The time derivative  of $V$ along the P-PDGD \eqref{eq:ppdgd} is 
\begin{align}
    \dot{V}(\bz) = \nabla_{\bz} V(\bz)^\top \dot{\bz} =   (\bz-\bz^* +H(\bz))^\top \bm{f}(\bz).
\end{align}
One useful property is stated as  Lemma \ref{lemma:1}.
\begin{lemma}\label{lemma:1}
For any $\alpha>0$, we have
\begin{align*}
    (\bz\!-\!\bz^* +\alpha H(\bz))^\top \bm{f}(\bz) \leq \!-||\bm{f}(\bz)||^2 \!- \alpha (\bz\!-\!\bz^*)^\top H(\bz).
\end{align*}
\end{lemma}
%
The proof of Lemma \ref{lemma:1} follows \cite[Lemma 2.4]{bansode2019exponential}. For completeness, here we provide a detailed proof  as 
the  three steps below:

\vspace{0.1cm}
\noindent 
1) We use the fact \cite{nagurney2012projected} that the projection operator satisfies 
\begin{align} \label{eq:step1}
     ( \mathrm{Proj}_{\sZ}(\bm{\gamma})-\bm{\beta})^\top (\bm{\gamma} -  \mathrm{Proj}_{\sZ}(\bm{\gamma}))\geq 0,
\end{align}
for all $ \bm{\gamma}\in\R^{|\sZ|}, \bm{\beta}\in \sZ$.

\noindent 
2) Let $\bm{\gamma} = \bz - \alpha H(\bz)$ and $\bm{\beta} = \bz^*$, then \eqref{eq:step1} becomes
\begin{align} \label{eq:step2}
     (\bm{f}({\bz})+\bz-\bz^*)^\top (\alpha H(\bz) + \bm{f}({\bz}) )\leq 0.
\end{align}

\noindent 
3) Thus we obtain Lemma \ref{lemma:1} by
\begin{align*}
  &  \, (\bz-\bz^* +\alpha H(\bz))^\top \bm{f}({\bz}) \\
  = &   \, (-\bm{f}({\bz}) +\bm{f}({\bz})+\bz -\bz^* +\alpha H(\bz))^\top \bm{f}({\bz}) \\
      = &  \, -\!||\bm{f}({\bz})||^2  + (\bm{f}({\bz})+\bz -\bz^*)^\top \bm{f}({\bz})+ \alpha H(\bz)^\top \bm{f}({\bz})\\
      \leq &   \, -\!||\bm{f}({\bz})||^2  - \alpha(\bm{f}({\bz})+\bz -\bz^*)^\top  H(\bz) + \alpha H(\bz)^\top \bm{f}({\bz})\\
      =&   \, -\!||\bm{f}({\bz})||^2  -  \alpha (\bz-{\bz}^*)^\top H(\bz)
\end{align*}
where the inequality above is because of (\ref{eq:step2}).

Using the result of Lemma \ref{lemma:1} with $\alpha=1$, we obtain 
\begin{align} \label{eq:lastineq}
    & \dot{V}(\bz)  \leq  - ||\bm{f}({\bz})||^2 - (\bz-\bz^*)^\top H(\bz) \nonumber\\
    = &\! -\! ||\bm{f}({\bz})||^2 \!-\!(\bx\!-\!\bx^*)^\top\nabla_{\bx}L(\bx,\!\bla)\!+\!(\bla\!-\!\bla^*)^\top\nabla_{\bla}L(\bx,\!\bla)  \nonumber \\
      \leq &\! -\! ||\bm{f}({\bz})||^2 +\! L(\bx^*,\bla) \!-\!L(\bx,\bla)  \!+\! L(\bx,\bla)\!-\!L(\bx,\bla^*) \nonumber   \\
      = &\!-\! ||\bm{f}({\bz})||^2 +\! L(\bx^*\!,\bla) \!-\!L(\bx^*\!,\bla^*) \! +\! L(\bx^*\!,\bla^*)\!-\!L(\bx,\bla^*) \nonumber \\
     \leq &\!-\! ||\bm{f}({\bz})||^2\leq 0
\end{align}
where the second inequality follows that $L(\bx,\bla)$ is convex in $\bx$ and concave in $\bla$. 

By \eqref{eq:lastineq}, we have that every compact level set of $V$ is forward invariant, and since $V$ is radially unbounded, it follows that all trajectories $\bz(t)$ remain bounded. Thus, by LaSalle's Theorem \cite[Theorem 4.4]{khalil2002nonlinear},  $\bz(t)$ converges to the largest invariant  compact subset $\mathcal{M}$ contained in $\mathcal{S}$:
\begin{align}
    \mathcal{S}:= \big\{{\bz}\in \sZ:~ \dot{V}({\bz})=0, \, V({\bz})\leq V(\bz(0))      \big\}. 
\end{align}
%
When $ \dot{V}({\bz})=0$,  we must have  $L(\bx^*,\bla) = L(\bx^*,\bla^*)$ and $L(\bx,\bla^*)=L(\bx^*,\bla^*)$ by (\ref{eq:lastineq}). Thus any point $\bz\in\mathcal{M}$ is an optimal solution of the saddle point problem (\ref{eq:saddle}).
 Lastly, the trick used in the proof of \cite[Theorem 15]{li2015connecting} can be adopted to show that $\bz(t)$ eventually converges to a fixed optimal point $\bz^*$. 
By strong duality (Assumption \ref{ass:finite}),
the component $\bx^*$ of the optimal point $\bz^*$ is the optimal solution of the OVC problem (\ref{eq:ovc}).
Thus Theorem \ref{thm:ppdgdcon} is proved.

\subsection{Proof of Theorem \ref{thm:spas}} \label{sec:app2}
Denote $\bs_1:=[\bx;\bla]$, $\bs_2:=[\bxi;\bmu]$, and $\bs:=[\bs_1;\bs_2]$. 
The ES-P-PDGD (\ref{eq:es:ppdgd}) is reformulated in compact form as
\begin{align} \label{eq:original}
 \dot{\bs}=    \begin{bmatrix}
      \dot{\bs}_1 \\   \dot{\bs}_2
    \end{bmatrix} = \begin{bmatrix}
    \bm{g}_1(\bs_1,\bs_2)\\
    \frac{1}{\epsilon}(-\bs_2 + \bm{g}_2(t, \bs_1))
    \end{bmatrix}:= \bm{g}(t, \bs),
\end{align}
where the function $\bg_1(\bs_1,\bs_2)$ captures the dynamics \eqref{eq:es:x}-\eqref{eq:es:la-}, and function $\bg_2(t,\bs_1)$ is given by
\begin{align}\label{eq:g2}
    \bg_2:=\begin{bmatrix}
      \big(  \frac{2}{a}  v_j(\bx\!+\! a\sin(\bom t )) \sin(\omega_n t )    \big)_{n\in [N],j\in\sM} \\
    \big(  v_j(\bx+ a\sin(\bom t)) \big)_{j\in\sM}
    \end{bmatrix},
\end{align}
where the first part and the second part are associated with the dynamics \eqref{eq:es:xi} of $\bxi$ and \eqref{eq:es:mu} of $\bmu$, respectively.

The following Lemma \ref{lemma:aveg2} states the average map for the function $\bg_2(t,\bs_1)$, which is proved in Appendix \ref{sec:app:lem}.

\begin{lemma}\label{lemma:aveg2}
The average of function $\bg_2(t,\bs_1)$ is given by
\begin{align}
    \bg_2^{\mathrm{av}}(\bs_1):=& \frac{1}{T}\!\int_0^T \!\bg_2(t,\bs_1)\, dt
    = \bm{\ell}(\bs_1) +\sO(a),
\end{align}
where 
$\bm{\ell}(\bs_1)\!:=\!\begin{bmatrix}(\frac{\partial v_j(\bx)}{\partial x_n})_{n\in[N],j\in\sM} \\
    (v_j(\bx))_{j\in\sM}
    \end{bmatrix}$, and $T$ is the  minimum common  period of the sinusoidal signals $\sin(\bom t)$.
\end{lemma}

We analyze the stability of the system \eqref{eq:original} via averaging theory and singular perturbation theory,
which is divided into the following three steps.

\vspace{0.1cm}
\noindent \textbf{Step 1)} \emph{Construct a compact  set to study the behavior of system \eqref{eq:original} restricted to it.}

To apply averaging theory and singular perturbation theory, it requires that the considered trajectories  stay within  predefined compact sets. Without loss of generality, we consider the compact set $[(\hat{\sA}+\Delta \mathbb{B})\cap \hat{\sZ}]\times \Delta \B$ for the initial condition $\bs(0)$ and any desired $\Delta>0$
. Here, 
$\B$ denotes a closed unit ball of appropriate dimension, and $\hat{\sA}+\Delta \mathbb{B}$ denotes 
 the union of all sets obtained by taking a closed ball of radius $\Delta$ around each point in the saddle point set $\hat{\sA}$. 

According to Theorem \ref{thm:ppdgdcon}, there exists a class-$\mathcal{KL}$ function $\beta$ such that for any initial condition $\bz(0)\in \hat{\sZ}$, the trajectory $\bz(t)$ of the P-PDGD \eqref{eq:ppdgd} with the feasible set $\hat{\sX}$ satisfies 
\begin{align} \label{eq:betanew}
     ||\bz(t)||_{\hat{\sA}}\leq \beta(||\bz(0)||_{\hat{\sA}},\,  t) , \quad \forall t\geq 0.
\end{align}

Without loss of generality, we assume the desired convergence precision $\nu\in(0,1)$. Using the $\beta$ function in \eqref{eq:betanew}, we define the set 
\begin{align}\label{eq:Fset}
    \sF\!:=\! \Big\{\bs_1\!\in \!\hat{\sZ}:   ||\bs_1||_{\hat{\sA}} \leq \beta\big( \max_{\bp\in \hat{\sA}+\Delta \mathbb{B} } ||\bp||_{\hat{\sA}},0 \big)+1  \Big\},\!
\end{align}
which is compact. Due to the boundedness of $\sF$, there exists a positive constant $M_1$ such that $\sF\subset M_1 \B$. Since $\bm{\ell}(\bs_1)$ (defined in Lemma \ref{lemma:aveg2}) is continuous by Assumption \ref{ass:con_sm}, 
 there exists a positive constant $M_2>\Delta$ such that $||\bm{\ell}(\bs_1)|| +1\leq M_2$ whenever $||\bs_1||\leq M_1$.   
 We then study the behavior of system \eqref{eq:original} \textbf{restricted to evolve in the compact set $\sF\times M_2\B$}.  

\vspace{0.1cm}
\noindent \textbf{Step 2)} \emph{Study the stability properties of the average system of the original system \eqref{eq:original}.}

By definition \eqref{eq:fredef}, the sinusoidal signals in system (\ref{eq:original}) are given by $\sin(\frac{2\pi}{\varepsilon_\omega} \kappa_n t )$ for $n\in[N]$. For sufficiently small $\varepsilon_\omega$, system (\ref{eq:original}), evolving in $\sF\times M_2\B$, is in  standard form  for the application of averaging theory. 
%
%
%
%
By Lemma \ref{lemma:aveg2},
we  derive the autonomous \textbf{average system}  of 
system \eqref{eq:original}, which is given by \eqref{eq:realave} (evolving in $\sF\times M_2\B$):
\begin{align} \label{eq:realave}
    \dot{\by} = \begin{bmatrix} 
    \dot{\by}_1\\
    \dot{\by}_2
    \end{bmatrix} = \frac{1}{T}\!\int_0^T \!\!\bg(t,\by)\, dt
    =\begin{bmatrix}
    \bm{g}_1(\by_1,\by_2)\\
    \frac{1}{\epsilon}(-\by_2 \!+\! \bm{\ell}(\by_1) \!+\!\sO(a)  )
    \end{bmatrix}
\end{align}
where $\by:=[\by_1;\by_2]$ takes the same form as $\bs:=[\bs_1;\bs_2]$.

To analyze the average system \eqref{eq:realave}, we can first ignore the small  $\sO(a)$-perturbation by setting  $a=0$. Thus the resultant system is in the standard form for the application of singular perturbation theory \cite{teel2003unified} with the slow dynamics of $\by_1$ and fast dynamics of $\by_2$.
As $\epsilon\to 0^+$, we freeze the slow state $\by_1$, and
the \textbf{boundary layer system} of the average system \eqref{eq:realave} with $a=0$ in the time scale $\tau=t/\epsilon$ is
\begin{align}
    \frac{d \by_2}{d \tau} = -\by_2 +\bm{\ell}(\by_1),
\end{align}
which is a linear time-invariant   system with the unique equilibrium point $\by_2^*=\bm{\ell}(\by_1)$. As a result,
the associated \textbf{reduced system}  is derived as 
\begin{align}\label{eq:reduced}
    \dot{\by}_1 = \bg_1(\by_1, \bm{\ell}(\by_1)),
\end{align}
which is precisely the P-PDGD \eqref{eq:ppdgd}. By Theorem \ref{thm:ppdgdcon} and \cite[Theorem 2]{Wang:12_Automatica}, it follows that as $\epsilon\to 0^+$, the set $\hat{\sA}\times {M}_2\mathbb{B}$ is semi-globally practically asymptotically stable (\textbf{SGPAS}) for the average system \eqref{eq:realave} with  $a=0$.
Then
by the structural robustness results
for ordinary differential equations with continuous right-hand sides \cite[Proposition A.1]{PovedaNaLi2019},  the set $\hat{\sA}\times {M}_2\mathbb{B}$ is also SGPAS for the average system \eqref{eq:realave}  as $(\epsilon,a)\to 0^+$, which is stated as Lemma \ref{lemma:y:sgpsa}.
\begin{lemma}\label{lemma:y:sgpsa} Given the precision $\nu$,
there exists $\epsilon^*>0$ such that for any $\epsilon\in(0,\epsilon^*)$, there exists $a^*>0$ such that for any $ a\in (0,a^*)$, with initial condition $\by(0)\in [(\hat{\sA}+\Delta \mathbb{B})\cap \hat{\sZ}]\times \Delta \B$,
the solution $\by(t)$ of the average system \eqref{eq:realave} satisfies that for all $ t\geq 0$,
\begin{align} \label{eq:y1con}
     ||\by_1(t)||_{\hat{\sA}}\leq \beta(||\by_1(0)||_{\hat{\sA}},\,  t) +\frac{\nu}{2}.
\end{align}
\end{lemma}
\noindent The proof of Lemma \ref{lemma:y:sgpsa} is provided in Appendix \ref{sec:app:lem:sgpas}.


\vspace{0.1cm}
\noindent \textbf{Step 3)} \emph{Link the stability property of the average system \eqref{eq:realave} to the stability property of the original system \eqref{eq:original}.}

Since  the set $\hat{A}\times M_2\B$ is SGAPAS for the average system \eqref{eq:realave}  as $(\epsilon,a)\to 0^+$, by averaging theory for perturbed systems \cite[Theorem 7]{PovedaNaLi2019}, it directly obtains that for each pair of $(\epsilon,a)$ inducing the bound \eqref{eq:y1con}, there exists $\varepsilon_\omega^*>0$ such that for any $\varepsilon_\omega\in(0,\varepsilon_\omega^*)$, the solution $\bs(t)$ of the original system \eqref{eq:original} restricted to $\sF\times M_2\B$ satisfies 
\begin{align}
    ||\bs_1(t)||_{\hat{\sA}}   \leq \beta(||\bs_1(0)||_{\hat{\sA}},\,  t) +\nu, \ \ \forall t\geq 0.
\end{align}
The completeness of solution $\bs$ for the original system \eqref{eq:original} is guaranteed by taking $M_2$ sufficiently large.

Thus Theorem \ref{thm:spas} is proved.

\subsection{Proof of Lemma \ref{lemma:aveg2}} \label{sec:app:lem}

We first consider the  integration on the first part of $\bg_2(t,\bs_1)$. By the Taylor expansion of $v_j(\cdot)$, each component of this integration is  ($\forall j\in\sM, n\in[N]$)
\begin{align*}
 & \frac{1}{T}\int_{0}^T \frac{2}{a}  v_j(\bx+ a\sin(\bom t )) \sin(\omega_n t )\, dt \\
 = &  \frac{1}{T}\int_{0}^T 
\frac{2}{a} \big[v_j(\bx)\!+\! a  \nabla v_j(\bx)^\top \sin (\bom t)
\! +\!  \sO(a^2)\big] \sin(\omega_n t )\, dt \\
= & \frac{1}{T}\int_{0}^T 2 \sum_{i=1}^N \frac{\partial v_j(\bx)}{\partial x_i} \sin(\omega_i t)\sin(\omega_n t )\, dt \!+\!\sO(a)\\
= & \frac{\partial v_j(\bx)}{\partial x_n} \frac{1}{T}\int_{0}^T \!  2\sin(\omega_n t )^2 \, dt +\sO(a)= \frac{\partial v_j(\bx)}{\partial x_n}  +\sO(a).
\end{align*}

As for the integration on the second part of $\bg_2(t,\bs_1)$, similarly, each component of this integration is ($\forall j\in\sM$)
\begin{align*}
     &\, \frac{1}{T}\int_{0}^T \!v_j(\bx+ a\sin(\bom t)) \, dt \\
= & \,\frac{1}{T}\int_{0}^T \!
v_j(\bx) + a  \nabla v_j(\bx)^\top \sin (\bom t) +\sO(a^2)\, dt\\
= & \,v_j(\bx) + \sO(a^2). 
\end{align*}
Combining these two parts, Lemma \ref{lemma:aveg2} is proved.

\subsection{Proof of Lemma \ref{lemma:y:sgpsa}} \label{sec:app:lem:sgpas}

By  the arguments of singular perturbation  and structural robustness (right above Lemma \ref{lemma:y:sgpsa}), 
it  follows that the bound \eqref{eq:y1con} holds for all $t\in [0, T_{\by})$,
where $[0, T_{\by})$ denotes the maximal time interval of existence of solution $\by$.

We further show that the solution $\by$ of the average system \eqref{eq:realave} exists for an unbounded time domain by the following lemma \ref{lemma:forward}, which follows a special case of \cite[Lemma 5]{NunoShamma_2020}.
\begin{lemma} \label{lemma:forward}
Let $M_2>0$ be given and $\bm{e}:\R_+ \to M_2\B$. Then for any $k>0$, the set $M_2\B$ is forward invariant for the dynamics $\dot{\bs}_2 = k(-\bs_2 +\bm{e}(t))$.
\end{lemma}

By the construction of $\sF$ \eqref{eq:Fset}, we obtain $\by_1(t)\in \mathrm{int}(\sF)$ for all $t\in [0, T_{\by})$.  Moreover, 
by setting $a^*$ sufficiently small such that
$||\sO(a)||<1$ for any $a\in(0,a^*)$, 
it follows that $||\by_1(t)||\leq M_1$ and $||\bm{\ell}(\by_1(t)) +\sO(a)|| < M_2$  for all $t\in[0, T_{\by})$. By Lemma \ref{lemma:forward}, it implies that $\by_2(t)\in \mathrm{int}(M_2\B)$ for all $t\geq 0$. Hence, the solution $\by(t)\in \mathrm{int}(\sF\times M_2 \B)$ for all $t\geq 0$, and has an unbounded time domain, i.e., $T_{\by}\to +\infty$.

\bibliography{IEEEabrv, ref}

\begin{thebibliography}{10}
\providecommand{\url}[1]{#1}
\csname url@samestyle\endcsname
\providecommand{\newblock}{\relax}
\providecommand{\bibinfo}[2]{#2}
\providecommand{\BIBentrySTDinterwordspacing}{\spaceskip=0pt\relax}
\providecommand{\BIBentryALTinterwordstretchfactor}{4}
\providecommand{\BIBentryALTinterwordspacing}{\spaceskip=\fontdimen2\font plus
\BIBentryALTinterwordstretchfactor\fontdimen3\font minus
  \fontdimen4\font\relax}
\providecommand{\BIBforeignlanguage}[2]{{%
\expandafter\ifx\csname l@#1\endcsname\relax
\typeout{** WARNING: IEEEtran.bst: No hyphenation pattern has been}%
\typeout{** loaded for the language `#1'. Using the pattern for}%
\typeout{** the default language instead.}%
\else
\language=\csname l@#1\endcsname
\fi
#2}}
\providecommand{\BIBdecl}{\relax}
\BIBdecl

\bibitem{8636257}
H.~{Sun}, Q.~{Guo}, J.~{Qi}, and et~al., ``Review of challenges and research
  opportunities for voltage control in smart grids,'' \emph{IEEE Trans. Power
  Syst.}, vol.~34, no.~4, pp. 2790--2801, 2019.

\bibitem{zheng2017robust}
W.~Zheng, W.~Wu, B.~Zhang, and Y.~Wang, ``Robust reactive power optimisation
  and voltage control method for active distribution networks via dual
  time-scale coordination,'' \emph{IET Generation, Transmission \&
  Distribution}, vol.~11, no.~6, pp. 1461--1471, 2017.

\bibitem{7039295}
B.~A. {Robbins}, H.~{Zhu}, and A.~D. {Domínguez-García}, ``Optimal tap
  setting of voltage regulation transformers in unbalanced distribution
  systems,'' \emph{IEEE Trans. Power Syst.}, vol.~31, no.~1, pp. 256--267,
  2016.

\bibitem{7244261}
B.~A. {Robbins} and A.~D. {Domínguez-García}, ``Optimal reactive power
  dispatch for voltage regulation in unbalanced distribution systems,''
  \emph{IEEE Trans. Power Syst.}, vol.~31, no.~4, pp. 2903--2913, 2016.

\bibitem{8779692}
G.~{Qu} and N.~{Li}, ``Optimal distributed feedback voltage control under
  limited reactive power,'' \emph{IEEE Trans. Power Syst.}, vol.~35, no.~1, pp.
  315--331, 2020.

\bibitem{8268542}
H.~J. {Liu}, W.~{Shi}, and H.~{Zhu}, ``Hybrid voltage control in distribution
  networks under limited communication rates,'' \emph{IEEE Trans. Smart Grid},
  vol.~10, no.~3, pp. 2416--2427, 2019.

\bibitem{7361761}
H.~{Zhu} and H.~J. {Liu}, ``Fast local voltage control under limited reactive
  power: Optimality and stability analysis,'' \emph{IEEE Trans. Power Syst.},
  vol.~31, no.~5, pp. 3794--3803, 2016.

\bibitem{7028508}
N.~{Li}, G.~{Qu}, and M.~{Dahleh}, ``Real-time decentralized voltage control in
  distribution networks,'' in \emph{52nd Annual Allerton Confer. on Commun.,
  Control, and Computing (Allerton)}, 2014, pp. 582--588.

\bibitem{8412143}
J.~{Zhang}, Z.~{Chen}, C.~{He}, Z.~{Jiang}, and L.~{Guan}, ``Data-driven-based
  optimization for power system var-voltage sequential control,'' \emph{IEEE
  Trans. Ind. Informat.}, vol.~15, no.~4, pp. 2136--2145, 2019.

\bibitem{7348696}
H.~{Zhang}, J.~{Zhou}, Q.~{Sun}, J.~M. {Guerrero}, and D.~{Ma}, ``Data-driven
  control for interlinked ac/dc microgrids via model-free adaptive control and
  dual-droop control,'' \emph{IEEE Trans. Smart Grid}, vol.~8, no.~2, pp.
  557--571, 2017.

\bibitem{7540852}
C.~{Mugnier}, K.~{Christakou}, J.~{Jaton}, M.~{De Vivo}, M.~{Carpita}, and
  M.~{Paolone}, ``Model-less/measurement-based computation of voltage
  sensitivities in unbalanced electrical distribution networks,'' in \emph{2016
  Power Systems Computation Conference (PSCC)}, 2016, pp. 1--7.

\bibitem{8873667}
H.~{Xu}, A.~D. {Domínguez-García}, V.~V. {Veeravalli}, and P.~W. {Sauer},
  ``Data-driven voltage regulation in radial power distribution systems,''
  \emph{IEEE Trans. Power Syst.}, vol.~35, no.~3, pp. 2133--2143, 2020.

\bibitem{8944292}
W.~{Wang}, N.~{Yu}, Y.~{Gao}, and J.~{Shi}, ``Safe off-policy deep
  reinforcement learning algorithm for volt-var control in power distribution
  systems,'' \emph{IEEE Trans. Smart Grid}, vol.~11, no.~4, pp. 3008--3018,
  2020.

\bibitem{9143169}
Y.~{Zhang}, X.~{Wang}, J.~{Wang}, and Y.~{Zhang}, ``Deep reinforcement learning
  based volt-var optimization in smart distribution systems,'' \emph{IEEE
  Trans. Smart Grid}, vol.~12, no.~1, pp. 361--371, 2021.

\bibitem{9076841}
S.~{Wang}, J.~{Duan}, D.~{Shi}, C.~{Xu}, H.~{Li}, R.~{Diao}, and Z.~{Wang}, ``A
  data-driven multi-agent autonomous voltage control framework using deep
  reinforcement learning,'' \emph{IEEE Trans. Power Syst.}, vol.~35, no.~6, pp.
  4644--4654, 2020.

\bibitem{9274529}
H.~{Liu} and W.~{Wu}, ``Two-stage deep reinforcement learning for
  inverter-based volt-var control in active distribution networks,'' \emph{IEEE
  Trans. Smart Grid}, pp. 1--1, 2020.

\bibitem{chen2021reinforcement}
X.~Chen, G.~Qu, Y.~Tang, S.~Low, and N.~Li, ``Reinforcement learning for
  decision-making and control in power systems: Tutorial, review, and vision,''
  \emph{arXiv preprint arXiv:2102.01168}, 2021.

\bibitem{9147814}
Y.~{Chen}, A.~{Bernstein}, A.~{Devraj}, and S.~{Meyn}, ``Model-free primal-dual
  methods for network optimization with application to real-time optimal power
  flow,'' in \emph{2020 American Control Conference (ACC)}, 2020, pp.
  3140--3147.

\bibitem{ariyur2003real}
K.~B. Ariyur and M.~Krsti{\'c}, \emph{Real Time Optimization by Extremum
  Seeking Control}.\hskip 1em plus 0.5em minus 0.4em\relax Wiley Online
  Library, 2003.

\bibitem{7397989}
M.~{Ye} and G.~{Hu}, ``Distributed extremum seeking for constrained networked
  optimization and its application to energy consumption control in smart
  grid,'' \emph{IEEE Trans. Control Syst. Technol.}, vol.~24, no.~6, pp.
  2048--2058, 2016.

\bibitem{8888186}
M.~D. {Sankur}, R.~{Dobbe}, A.~{von Meier}, and D.~B. {Arnold}, ``Model-free
  optimal voltage phasor regulation in unbalanced distribution systems,''
  \emph{IEEE Trans. Smart Grid}, vol.~11, no.~1, pp. 884--894, 2020.

\bibitem{6362193}
X.~{Li}, Y.~{Li}, and J.~E. {Seem}, ``Maximum power point tracking for
  photovoltaic system using adaptive extremum seeking control,'' \emph{IEEE
  Trans. Control Syst. Technol.}, vol.~21, no.~6, pp. 2315--2322, 2013.

\bibitem{7350258}
D.~B. {Arnold}, M.~{Negrete-Pincetic}, M.~D. {Sankur}, D.~M. {Auslander}, and
  D.~S. {Callaway}, ``Model-free optimal control of var resources in
  distribution systems: An extremum seeking approach,'' \emph{IEEE Trans. Power
  Syst.}, vol.~31, no.~5, pp. 3583--3593, 2016.

\bibitem{8772122}
H.~{Nazaripouya}, H.~R. {Pota}, C.~{Chu}, and R.~{Gadh}, ``Real-time model-free
  coordination of active and reactive powers of distributed energy resources to
  improve voltage regulation in distribution systems,'' \emph{IEEE Trans.
  Sustain. Energy}, vol.~11, no.~3, pp. 1483--1494, 2020.

\bibitem{8478426}
J.~{Johnson}, A.~{Summers}, R.~{Darbali-Zamora}, J.~{Hernandez-Alvidrez},
  J.~{Quiroz}, D.~{Arnold}, and J.~{Anandan}, ``Distribution voltage regulation
  using extremum seeking control with power hardware-in-the-loop,'' \emph{IEEE
  J. Photovolt.}, vol.~8, no.~6, pp. 1824--1832, 2018.

\bibitem{flaxman2004online}
A.~D. Flaxman, A.~T. Kalai, and H.~B. McMahan, ``Online convex optimization in
  the bandit setting: gradient descent without a gradient,'' \emph{arXiv
  preprint cs/0408007}, 2004.

\bibitem{Manzie09}
C.~Manzie and M.~Krstic, ``Extremum seeking with stochastic perturbations,''
  \emph{IEEE Trans. Autom. Control}, vol.~54, no.~3, pp. 580--585, 2009.

\bibitem{Teel01}
A.~R. Teel and D.~Popovic, ``Solving smooth and nonsmooth multivariable
  extremum seeking problems by the methods of nonlinear programming,'' \emph{In
  Proc. of American Control Conference}, pp. 2394--2399, 2001.

\bibitem{Khong13}
S.~Khong, D.~Nesic, Y.~Tan, and C.~Manzie, ``Unified frameworks for
  sampled-data extremum seeking control: Global optimisation and multi-unit
  systems,'' \emph{Automatica}, no.~49, pp. 2720--2733, 2013.

\bibitem{PovedaTAC17}
J.~I. Poveda and A.~R. Teel, ``A robust event-triggered approach for fast
  sampled-data extremization and learning,'' \emph{IEEE Trans. Autom. Control},
  no.~10, pp. 4949--4964, 2017.

\bibitem{gao2003exponential}
X.-B. Gao, ``Exponential stability of globally projected dynamic systems,''
  \emph{IEEE Trans. Neural Netw.}, vol.~14, no.~2, pp. 426--431, 2003.

\bibitem{nagurney2012projected}
A.~Nagurney and D.~Zhang, \emph{Projected Dynamical Systems and Variational
  Inequalities with Applications}.\hskip 1em plus 0.5em minus 0.4em\relax
  Springer Science \& Business Media, 2012, vol.~2.

\bibitem{8571158}
Y.~{Zhu}, W.~{Yu}, G.~{Wen}, and G.~{Chen}, ``Projected primal–dual dynamics
  for distributed constrained nonsmooth convex optimization,'' \emph{IEEE
  Trans. Cybern.}, vol.~50, no.~4, pp. 1776--1782, 2020.

\bibitem{bansode2019exponential}
P.~Bansode, V.~Chinde, S.~Wagh, R.~Pasumarthy, and N.~Singh, ``On the
  exponential stability of projected primal-dual dynamics on a riemannian
  manifold,'' \emph{arXiv preprint arXiv:1905.04521}, 2019.

\bibitem{khalil2002nonlinear}
H.~K. Khalil and J.~W. Grizzle, \emph{Nonlinear Systems}, 3rd~ed.\hskip 1em
  plus 0.5em minus 0.4em\relax Prentice hall Upper Saddle River, NJ, 2002.

\bibitem{Wang:12_Automatica}
W.~Wang, A.~Teel, and D.~Ne\u{s}i\'{c}, ``Analysis for a class of singularly
  perturbed hybrid systems via averaging,'' \emph{Automatica}, vol.~48, no.~6,
  pp. 1057--1068, 2012.

\bibitem{qu2018exponential}
G.~Qu and N.~Li, ``On the exponential stability of primal-dual gradient
  dynamics,'' \emph{IEEE Contr. Syst. Lett.}, vol.~3, no.~1, pp. 43--48, 2018.

\bibitem{6760555}
M.~{Farivar}, L.~{Chen}, and S.~{Low}, ``Equilibrium and dynamics of local
  voltage control in distribution systems,'' in \emph{52nd IEEE Conference on
  Decision and Control}, 2013, pp. 4329--4334.

\bibitem{PovedaNaLi2019}
J.~I. Poveda and N.~Li, ``Robust hybrid zero-order optimization algorithms with
  acceleration via averaging in time,'' \emph{Automatica}, vol. 123, p. 109361,
  2021.

\bibitem{matpower}
R.~D. Zimmerman, C.~E. Murillo-S{\'a}nchez, and R.~J. Thomas, ``Matpower:
  Steady-state operations, planning and analysis tools for power systems
  research and education,'' \emph{IEEE Trans. Power Syst.}, vol.~26, no.~1, pp.
  12--19, Feb. 2011.

\bibitem{cvx}
M.~Grant and S.~Boyd, ``{CVX}: Matlab software for disciplined convex
  programming, version 2.1,'' \url{http://cvxr.com/cvx}, Mar. 2014.

\bibitem{ruszczynski2011nonlinear}
A.~Ruszczynski, \emph{Nonlinear Optimization}.\hskip 1em plus 0.5em minus
  0.4em\relax Princeton university press, 2011.

\bibitem{clarke1990optimization}
F.~H. Clarke, \emph{Optimization and Nonsmooth Analysis}.\hskip 1em plus 0.5em
  minus 0.4em\relax Wiley: Society Series of Monographs and Advanced Texts,
  SIAM, 1990.

\bibitem{li2015connecting}
N.~Li, C.~Zhao, and L.~Chen, ``Connecting automatic generation control and
  economic dispatch from an optimization view,'' \emph{IEEE Control Netw.
  Syst.}, vol.~3, no.~3, pp. 254--264, 2015.

\bibitem{teel2003unified}
A.~R. Teel, L.~Moreau, and D.~Nesic, ``A unified framework for input-to-state
  stability in systems with two time scales,'' \emph{IEEE Trans. Autom.
  Control}, vol.~48, no.~9, pp. 1526--1544, 2003.

\bibitem{NunoShamma_2020}
S.~Park, N.~Martins, and J.~Shamma, ``Payoff dynamics model and evolutionary
  dynamics model: Feedback and convergence to equilibria,''
  \emph{arXiv:1903.02018v4}, 2020.

\end{thebibliography}

\end{document}